\documentclass[11pt]{article}
\usepackage[T1]{fontenc}
\usepackage{lmodern}
\usepackage[margin=1in]{geometry}
\usepackage{amsmath,amssymb,amsthm}
\usepackage{microtype,enumitem,graphicx}
\usepackage{tikz}
\usepackage[hidelinks]{hyperref}
\hypersetup{
pdftitle={Uniqueness of asymmetric incompressible jets},
pdfauthor={Jianfeng Cheng and Zikang Gu}
}
\numberwithin{equation}{section}
\newtheorem{theorem}{Theorem}[section]
\newtheorem{proposition}[theorem]{Proposition}
\newtheorem{lemma}[theorem]{Lemma}

\theoremstyle{definition}
\newtheorem{definition}[theorem]{Definition}
\theoremstyle{remark}
\newtheorem{remark}[theorem]{Remark}
\newcommand{\R}{\mathbb R}
\setlist[enumerate]{leftmargin=*,itemsep=5pt,topsep=6pt}
\allowdisplaybreaks[1]
\title{Complete uniqueness for two-dimensional asymmetric
		incompressible jet flows}
\author{%
  \normalsize
  \textsc{Jianfeng Cheng},
  \textsc{Zikang Gu}\textsuperscript{$\dagger$}\\[2.5em]
  {\normalsize Department of Mathematics, Sichuan University,}\\[-1pt]
  {\normalsize Chengdu 610065, P.~R.~China.}%
}
\date{}
\begin{document}
\maketitle

\begingroup
\renewcommand{\thefootnote}{\fnsymbol{footnote}}
\footnotetext[1]{%
  \emph{This work is supported in part by National Key R\&D Program of China (No. 2022YFA1007700), National Natural Science Foundation of China (No. 12371237 and No. 12526202) and Sichuan Provincial Natural Science Foundation (No. 2026NSFSCZY0123).}%
  \par\noindent
  E-Mail: \href{mailto:jfcheng@scu.edu.cn}{jfcheng@scu.edu.cn}
  (J.~Cheng), E-Mail:
  \href{mailto:zik2019@163.com}{zik2019@163.com} (Z.~Gu).%
  \par\noindent
  \textsuperscript{$\dagger$}Corresponding author: Zikang Gu.
  E-Mail: \href{mailto:zik2019@163.com}{zik2019@163.com}.%
}
\endgroup

\begin{abstract}
	Alt, Caffarelli and Friedman established the existence of two-dimensional asymmetric incompressible jet flows in [\emph{Asymmetric jet flows}. Comm. Pure Appl. Math. \textbf{35}, 29--68, 1982]. In their formulation, the constant fluid speed $\lambda$ on the free boundaries and the downstream asymptotic angle $\theta$ are selected as part of the solution. In this paper, we prove that the asymmetric jet is unique for a prescribed nozzle and mass flux. We first establish the strict vertical monotonicity of the stream function $u$ and introduce the partial hodograph coordinates $\xi=x$ and $\eta=u(x,y)$, which map every fluid domain onto the same fixed strip. Under this framework, each jet is characterized by a height function $y=h(\xi,\eta)$, which determines the vertical position along the streamline $u=\eta$ at horizontal coordinate $\xi$. By analyzing the difference between the height functions of two candidate jets on this strip, we derive an elliptic equation with Dirichlet conditions on the fixed walls and conormal conditions on the free boundaries. Combining these mixed boundary conditions with  the asymptotic behaviors, we first establish the uniqueness of the free-boundary fluid speed $\lambda$. For the resulting common speed, a local perturbation of the height function near a nozzle endpoint, together with the maximum principle, proves that the two height functions must coincide. Inverting the partial hodograph transformation then establishes the uniqueness of the physical jet.
\end{abstract}

\section{Introduction}

\subsection{Background and the uniqueness problem}
\label{sec-background}

Free-streamline theory has its origin in the classical study of
liquid jets, wakes, and cavities. Kirchhoff \cite{Kirchhoff1869}
gave an early analytical treatment of incompressible jet flows
through a nozzle, and systematic accounts of the subsequent
development were provided by Birkhoff and Zarantonello
\cite{BZ1957}, Gilbarg \cite{Gilbarg1960}, and Wu \cite{Wu1972}.
For a jet issuing from a nozzle, part of the fluid boundary is
prescribed by the solid walls, while the remaining part consists
of free streamlines that must be determined together with the
velocity field. Under the assumptions of two-dimensional,
irrotational, steady incompressible flow, the stream function
is harmonic. Moreover, if gravity and surface tension are
neglected and the exterior pressure is constant, Bernoulli's
law requires the flow speed to be constant on the free
boundaries. Thus, the simultaneous determination of the stream
function, the fluid domain, the free boundaries, and the flow
speed on these boundaries makes the problem nonlinear.

The early existence theory for incompressible jet flows was
developed by Weinstein \cite{Weinstein1930}, Leray
\cite{Leray1935}, and Finn \cite{Finn1956}. These works relied
mainly on hodograph transformations, conformal mappings,
fixed-point arguments, and approximation of curved walls by
polygonal ones. A different approach to free boundary problems
was introduced by Alt and Caffarelli \cite{AC81} through their
variational method. Alt, Caffarelli and Friedman subsequently
applied this approach to establish the existence of solutions
to incompressible free-streamline problems, including
two-dimensional asymmetric jets \cite{ACF}, jets under gravity
\cite{ACFgravity}, axially symmetric jets \cite{ACFaxis}, and
two-fluid jets \cite{ACFtwofluidsI,ACFtwofluidsII}. Their
variational method made it possible to study nozzle walls
that were not covered by the convexity assumptions of the
earlier theory.

The variational framework of Alt, Caffarelli and Friedman has
since been adapted to many nozzle configurations and fluid
models. In the two-dimensional symmetric setting, Cheng,
Du and Xiang \cite{CDX19} proved the existence and uniqueness
of incompressible jet flow through a de Laval nozzle with
a prescribed flow speed on the free boundary. They also
showed that the detachment point varies monotonically with
this speed. Du and Tang \cite{DuTang21} established the
existence of asymmetric incompressible jet flow through
a finitely long nozzle with prescribed inlet velocity and
surrounding pressure. Further well-posedness results were
obtained for jets flowing past an obstacle
\cite{ChengHuang23}, two-dimensional rotational jets with
gravity \cite{ChengPei24}, and axially symmetric rotational
jets \cite{CGX25}. For further results on jet flows, we
refer the reader to
\cite{ACFcompressible,Berg1962,CDXcompressible19,LSTX25,WangXin19}
and the references therein.

We focus here on the two-dimensional asymmetric incompressible
jet problem for which Alt, Caffarelli and Friedman
\cite[Theorem~1.1]{ACF} established the existence of solutions
by a variational approach. In their formulation, a solution
consists of a fluid domain $\Omega$, a harmonic stream
function $u$, two detached free streamlines $\Gamma_1$
and $\Gamma_2$, a free-boundary speed $\lambda>0$,
and a downstream asymptotic direction
$\mathbf{e}=(\cos\theta,\sin\theta)$ with
$\theta\in(-\pi/2,\pi/2)$, satisfying
\[
\Delta u=0\quad\text{in }\Omega,\qquad
u=1\quad\text{on }N_1\cup\Gamma_1,\qquad
u=-1\quad\text{on }N_2\cup\Gamma_2,
\qquad
|\nabla u|=\lambda\quad\text{on }\Gamma_1\cup\Gamma_2,
\]
where $N_1$ and $N_2$ denote the upper and lower solid walls,
respectively. As $x\to+\infty$, the free boundaries
$\Gamma_1$ and $\Gamma_2$ approach two parallel lines
oriented along $\mathbf{e}$. Since the normalized flux
is $2$, conservation of mass requires the distance between
these lines to be $2/\lambda$. Neither $\lambda$ nor
$\theta$ is prescribed a priori. Both parameters are
determined together with the free boundaries $\Gamma_1$
and $\Gamma_2$, which must detach smoothly from the
endpoints of $N_1$ and $N_2$, respectively.

As pointed out in \cite[p.~65]{ACF}, every minimizer
in their variational framework gives a solution of the
asymmetric jet problem, but whether every jet solution
arises as a variational minimizer remained open.
Thus, their variational construction did not settle
uniqueness among all solutions of the asymmetric jet problem. In fact, the lack of symmetry of the semi-infinitely long nozzle, coupled with the undetermined detachment speed $\lambda$ and asymptotic angle $\theta$, poses an obstacle to proving uniqueness of the solution of asymmetric jet flows.

The uniqueness theory for free-streamline problems goes
back to Weinstein's early work \cite{Weinstein1924}.
A systematic comparison method was subsequently developed
by Gilbarg \cite{Gilbarg1952} for axially symmetric jet
flows from nozzles that are star-shaped with respect
to a point on the axis
(see also \cite[Theorem~1.1]{Friedman1984}).
Serrin
\cite{Serrin1952Uniqueness,Serrin1952Comparison,Serrin1953}
extended this method to broader classes of plane and
axisymmetric incompressible free-streamline problems.
For planar symmetric nozzles, Finn \cite{Finn1956}
proved uniqueness under the assumption of polygonal
convex geometry. These classical comparison arguments rely on global
geometric assumptions, such as convexity or star-shapedness
of the nozzle walls. Similarly, the uniqueness results
for the related jet models discussed above use
two-dimensional or axial symmetry, involve only one
unknown free boundary, or prescribe additional inlet,
outlet, pressure, or detachment data. For a general
two-dimensional asymmetric nozzle considered in \cite{ACF}, however, the upper
and lower walls need not be related by reflection
symmetry. Hence, the two asymmetric free boundaries must be
determined simultaneously, together with the unknown
speed $\lambda$ and asymptotic angle $\theta$.
These features prevent a direct application of the
preceding uniqueness arguments.

A decisive step toward uniqueness of the solution to asymmetric jet problem was made by Cheng, Du, and Wang \cite{CDW}. For fixed free-boundary speed $\lambda$,
given two candidate solutions $u$ and $\widetilde u$ with asymptotic angles $\theta$ and $\widetilde\theta$, they proved that $u=\widetilde u$ and $\theta=\widetilde \theta$ by the sliding approach, provided that $\partial _yu\geq 0$ and $\partial_y \widetilde u\geq 0$.
However, the uniqueness of the free-boundary speed $\lambda$ for asymmetric jet flows remained an open problem \cite[Remark~1.1]{CDW}.
For asymmetric jet flows, neither the sliding method
in \cite{CDW} nor the scaling argument in \cite{Serrin1953}
can be directly applied to establish the uniqueness of
$\lambda$, since neither approach yields a contradiction
with the prescribed Neumann condition on the free boundaries.

In this paper, we prove the complete uniqueness of
two-dimensional asymmetric incompressible jet flows
for a prescribed nozzle and mass flux under the
assumptions stated in Subsection~\ref{sec-setting}.
We first establish the monotonicity of the stream function, which allows
us to use a partial hodograph transformation to
compare two candidate jets through their height
functions on a common fixed strip. Combining
comparison arguments with the far-field asymptotic
behavior, we prove that the free-boundary speeds
coincide. For this common speed, a local perturbation
of the height function near a nozzle endpoint,
together with the maximum principle and the
far-field asymptotic behavior, yields equality
of the height functions and hence of the physical
jets. The proof is outlined in
Subsection~\ref{sec-idea}.

\subsection{Mathematical formulation and the main result}
\label{sec-setting}

In this subsection, we formulate the geometric setup of the two-dimensional asymmetric nozzle, specify the governing equations for the stream function, define the notion of an asymmetric incompressible jet, and state the main uniqueness result.

\begin{figure}[htbp]
	\centering
	\includegraphics[width=80mm]{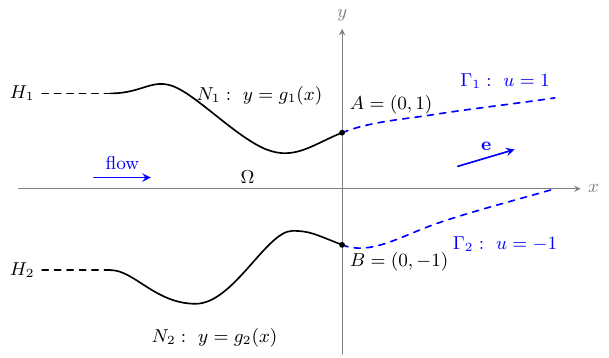}
	\caption{A two-dimensional asymmetric incompressible jet.}
	\label{fig-jet}
\end{figure}

\paragraph{The two-dimensional asymmetric nozzle.}
The upper and lower solid walls of the two-dimensional asymmetric nozzle are given respectively by the graphs
\[
N_1 = \{(x,g_1(x)) : x \le 0\} \quad \text{and} \quad N_2 = \{(x,g_2(x)) : x \le 0\}.
\]
They are assumed to satisfy the following geometric and regularity conditions:
\begin{enumerate}[label=(N\arabic*)]
	\item The continuous functions $g_1(x)$ and $g_2(x)$ are locally piecewise $C^{1,\beta}$ in $(-\infty,0)$ for some $\beta \in (0,1)$. Moreover,
	\[
	g_1(0) = 1, \quad g_2(0) = -1, \quad \text{and} \quad g_2(x) < g_1(x) \quad \text{for all } x \le 0.
	\]
	
	\item The profiles satisfy $g_1, g_2 \in C^{1,\beta}([-\delta,0])$ for some $\delta > 0$. In addition, there exist two asymptotic heights $H_2 < H_1$ such that
	\begin{equation}\label{nozzle-upstream}
		\lim_{x\to-\infty} g_i(x) = H_i,
		\quad \text{for } i = 1,2.
	\end{equation}
\end{enumerate}
We denote by $A = (0,1)$ and $B = (0,-1)$ the endpoints of the upper and lower nozzle walls, respectively. Figure~\ref{fig-jet} illustrates the nozzle geometry, the free boundaries, and the resulting jet flow region.

\paragraph{The stream function.}
We fix the Cartesian coordinate system, normalize the fluid density to unity, and normalize the stream function so that the total mass flux across any vertical cross-section equals two. Let $\mathbf{u} = (u_1, u_2)$ denote the normalized velocity field and let $p$ be the fluid pressure. The two-dimensional incompressible flow is governed by the following Euler equations
\begin{equation*}
	\left\{ \begin{array}{l}
		\partial_x u_1+\partial_y u_2=0,\\
		u_1\partial_x u_1+u_2\partial_y u_1+\partial_x p=0,\\
		u_1\partial_x u_2+u_2\partial_y u_2+\partial_y p=0.
	\end{array} \right.
\end{equation*}
The incompressibility condition allows us to introduce a
stream function $u(x,y)$ through the relations $u_1=u_y$ and
$u_2=-u_x$. Since the flow considered here is irrotational, we have
$\partial_x u_2-\partial_y u_1=0$, which gives
\[
\Delta u
=\partial_y u_1-\partial_x u_2
=0.
\]
Moreover, the definition of the stream function yields
\[
\mathbf{u}\cdot\nabla u
=u_yu_x-u_xu_y
=0.
\]
Thus, $u$ is constant along each streamline, and the
streamlines are given by the level curves of $u$. We now set the normalization $u = -1$ on the lower boundary $N_2\cup\Gamma_2$ and $u = 1$ on the upper boundary $N_1\cup\Gamma_1$. If the boundaries $N_2\cup\Gamma_2$ and $N_1\cup\Gamma_1$ are parameterized by $y = a(x)$ and $y = b(x)$ respectively, the total volume flux through any vertical cross-section is given by
\begin{equation}\label{physical-flux}
	\int_{a(x)}^{b(x)} u_1(x,y)\,dy = \int_{a(x)}^{b(x)} u_y(x,y)\,dy = u(x,b(x)) - u(x,a(x)) = 2.
\end{equation}
Neglecting gravity and surface tension, Bernoulli's law gives that
\[
p + \frac{1}{2}|\nabla u|^2 = \mathcal{B} \quad \text{in the fluid domain},
\]
for some constant $\mathcal{B} \in \mathbb{R}$. Under the condition of constant atmospheric pressure $p = p_{\mathrm{atm}}$ along the free streamlines, the flow speed must be constant on the free boundaries, namely,
\[
|\nabla u| = \lambda := \sqrt{2(\mathcal{B} - p_{\mathrm{atm}})} > 0.
\]
While the nozzle geometry and total mass flux are prescribed a priori, the free-boundary speed $\lambda$ is an unknown parameter to be determined as a part of the solution.

We now state the precise definition of the solution to two-dimensional asymmetric incompressible jet flows.

\begin{definition}[The solution to two-dimensional asymmetric incompressible jet flows]\label{def-jet}
	A solution of two-dimensional asymmetric incompressible jet flows is a quintuple $(u, \Gamma_1, \Gamma_2, \lambda,\theta)$ with $\theta\in(-\pi/2,\pi/2)$ satisfying the following conditions:
	\begin{enumerate}[label=(J\arabic*)]
		\item The upper and lower free boundaries are given respectively by $\Gamma_i = \{(x,f_i(x)) : x > 0\}$ for $i=1,2$, where
		\[
		f_i \in C^1([0,\infty)) \cap C^2((0,\infty)) \quad \text{for } i=1,2, \quad \text{and} \quad f_2(x) < f_1(x) \quad \text{for all } x \ge 0.
		\]
		Moreover, the continuous and smooth fit conditions at $A$ and $B$ hold, namely,
		\begin{equation}\label{jet-fit}
			g_i(0) = f_i(0), \qquad g_i'(0) = f_i'(0)=\lim_{x\rightarrow 0^+}f_i'(x) \quad \text{for } i=1,2.
		\end{equation}
		
		\item The stream function
        $u\in C(\overline\Omega)\cap
        C_{loc}^{1,\alpha}(\Omega\cup\Gamma_1\cup\Gamma_2)$
        for some $\alpha\in(0,1)$ satisfies
		\begin{equation}\label{jet-equations}
			\left\{
			\begin{aligned}
				&\Delta u = 0,\quad-1\leq u\leq 1 &&\text{in } \Omega, \\
				&u(x,a(x)) = -1, \quad u(x,b(x)) = 1 &&\text{for } x \in \mathbb{R}, \\
				&|\nabla u| = \lambda > 0 &&\text{on } \Gamma_1 \cup \Gamma_2,
			\end{aligned}
			\right.
		\end{equation}
		where the boundary functions $a(x)$ and $b(x)$ are defined by
		\begin{equation}\label{jet-ab}
			a(x) =
			\begin{cases}
				g_2(x), & x \le 0, \\
				f_2(x), & x > 0,
			\end{cases}
			\qquad
			b(x) =
			\begin{cases}
				g_1(x), & x \le 0, \\
				f_1(x), & x > 0,
			\end{cases}
		\end{equation}
		and the fluid domain is
		\begin{equation}\label{jet-domain}
			\Omega = \{(x,y) \in \mathbb{R}^2 : a(x) < y < b(x)\}.
		\end{equation}
	 Furthermore, there exists a small constant $r_0 > 0$ such that $u$ is $C^1$-smooth up to the endpoints $A$ and $B$:
		\begin{equation}\label{jet-endpoint-regularity}
			u \in C^1(\overline\Omega \cap B_{r_0}(A)) \cap C^1(\overline\Omega \cap B_{r_0}(B)).
		\end{equation}
		
		\item There exist constants  $b_1 > b_2$ such that the free boundaries satisfy the downstream asymptotic behaviors
		\begin{equation}\label{jet-downstream}
			\lim_{x\to+\infty} \bigl(f_i(x) - x\tan\theta - b_i\bigr) = 0, \qquad
			\lim_{x\to+\infty} f_i'(x) = \tan\theta \quad \text{for } i=1,2.
		\end{equation}
		Moreover, conservation of mass determines the limiting
		vertical distance between the two free boundaries
		in terms of $\lambda$, namely,
		\begin{equation}\label{jet-flux-width}
			\lambda(b_1 - b_2)\cos\theta=2.
		\end{equation}
	\end{enumerate}
\end{definition}

\begin{remark}
	\label{rem-normalization}
	In the formulation of Alt, Caffarelli, and Friedman \cite{ACF}, the unnormalized stream function $U$ takes the boundary values
	\[
	U = Q \quad \text{on } N_1 \cup \Gamma_1, \qquad U = -Q \quad \text{on } N_2 \cup \Gamma_2,
	\]
	for some prescribed constant $Q>0$, so that the total mass
	flux is $2Q$. The stream function in Definition~\ref{def-jet} corresponds to the normalization $u = U/Q$, which fixes the boundary data to $u = \pm 1$ and normalizes the mass flux to $2$. The unnormalized velocity field $(V_1, V_2)$ and the free-boundary flow speed $\Lambda$ are related to our normalized quantities by
	\[
	(V_1, V_2) = Q(u_y, -u_x), \qquad \Lambda = Q\lambda.
	\]
	This linear scaling is clearly invertible and preserves the underlying mathematical structure, while ensuring consistency with the notations in \cite{ACF}.
\end{remark}

\begin{remark}\label{rem-jet-definition}
	We briefly clarify the physical and geometric meaning of conditions \textup{(J1)--(J3)}:
	\begin{itemize}
		\item Condition \textup{(J1)} describes the detachment of the two unknown free boundaries and requires each streamline to detach tangentially and smoothly from the corresponding nozzle wall.
		\item In condition \textup{(J2)}, the constant Dirichlet values signify that the solid walls and free boundaries form impenetrable streamlines, with their difference equaling the normalized mass flux $2$ in view of \eqref{physical-flux}. The gradient constraint $|\nabla u| = \lambda$ reflects the constant-pressure condition through Bernoulli's law.
		\item Condition \textup{(J3)} characterizes the downstream asymptotic state: the flow asymptotically aligns with the unit direction $\mathbf{e} = (\cos\theta,\sin\theta)$, and the orthogonal width between the two asymptotic free boundaries is given by $(b_1-b_2)\cos\theta= 2/\lambda$.
	\end{itemize}
\end{remark}

The existence of the solution to two-dimensional asymmetric incompressible jet flows has been
established by Alt, Caffarelli and Friedman in \cite[Theorem 1.1]{ACF}.

{\bf Theorem A.} {\it  Assume that the nozzle walls satisfy conditions \textup{(N1)--(N2)}, then there exists a
solution $(u, \Gamma_1, \Gamma_2,\lambda,\theta)$ under the sense of definition \ref{def-jet} to two-dimensional asymmetric incompressible jet flows with $\lambda>0$ and $\theta\in\left(-\frac{\pi}2,\frac{\pi}2\right)$.}

\paragraph{Main result.}
Our main theorem establishes the complete uniqueness of asymmetric incompressible jets, asserting the coincidence of the free-boundary flow speeds, asymptotic downstream parameters, free boundaries, and the entire flow fields.

\begin{theorem}
	\label{thm-main}
	 Let $
	(u,\Gamma_1, \Gamma_2, \lambda,\theta)$ and $(\tilde u,\tilde \Gamma_1,\tilde \Gamma_2, \tilde\lambda,\tilde\theta)$
	be two-dimensional asymmetric incompressible jet flows. Then the two solutions are identical:
	\[
	u = \tilde u,\ \ \  \lambda = \tilde\lambda, \ \ \  \theta = \tilde\theta \ \ \ \text{and}\ \ \ \Gamma_i = \tilde \Gamma_i \quad \text{for } i=1,2.
	\]
\end{theorem}
\begin{remark}	
Gilbarg and Serrin
\cite{Gilbarg1952,Serrin1952Uniqueness,Serrin1952Comparison,Serrin1953}
studied the uniqueness of solutions to two-dimensional
and axially symmetric incompressible free-streamline problems.
They first established Under-Over theorem and comparison theorem of linear elliptic equation with one free boundary, and proved the single intersection property of the free boundary, namely, if a straight line does not intersect the nozzle wall, then it can intersect the free boundary at most once. This implies that the free boundary is star-shaped, provided that the symmetric nozzle wall is star-shaped with respect to a point on the axis.
The resulting geometric structure makes it possible to
compare two candidate solutions by scaling and thereby
establish uniqueness. However, for two-dimensional
asymmetric jet flows, the presence of two free boundaries
prevents a direct application of this scaling argument.
\end{remark}
\begin{remark}	
In the previous work \cite{CDW}, Cheng, Du, and Wang
established the uniqueness of the jet solution
$(u,\Gamma_1,\Gamma_2,\theta)$ for a prescribed
free-boundary speed $\lambda$. Their proof is based on
a sliding argument. To explain this approach, we consider
two candidate solutions with stream functions $u$ and
$\widetilde u$, free-boundary speeds $\lambda$ and
$\widetilde\lambda$, and asymptotic angles $\theta$
and $\widetilde\theta$, respectively. By translating
the first stream function vertically according to
$u^\varepsilon(x,y)=u(x,y-\varepsilon)$ for
$\varepsilon\geq0$, one can decrease $\varepsilon$
until the corresponding free boundaries first touch.
At a contact point, the Hopf boundary point lemma
yields a strict inequality between the normal
derivatives of $u^\varepsilon$ and $\widetilde u$.
When $\lambda=\widetilde\lambda$, this contradicts
the free-boundary conditions
$|\nabla u|=\lambda$ and
$|\nabla\widetilde u|=\widetilde\lambda$.
However, if $\lambda\ne\widetilde\lambda$,
the strict inequality between the normal derivatives
may be consistent with these two free-boundary
conditions. Therefore, the sliding argument
establishes uniqueness for a prescribed $\lambda$,
but does not exclude the existence of solutions
with different free-boundary speeds.
\end{remark}

\subsection{Strategy of the proof}
\label{sec-idea}

Because the stream functions of two candidate jets may be defined on different physical domains, their governing equations cannot be compared directly in physical space. To circumvent this difficulty, we transform the problem onto a fixed reference domain via a partial hodograph transformation.

The foundation of this transformation is the vertical monotonicity of the stream function $u$, that is, $u_y > 0$ throughout $\Omega$. We establish this by analyzing the vertical translation differences $u(x,y+t) - u(x,y)$ for small $t > 0$, applying a maximum principle on unbounded domains to deduce nonnegativity, and then invoking the strong maximum principle to ensure strict positivity. Consequently, for each fixed horizontal position $x = \xi$, the vertical profile $y \mapsto u(\xi,y)$ is a strictly increasing homeomorphism from the flow interval $[a(\xi), b(\xi)]$ onto $[-1, 1]$. We can therefore define the partial hodograph coordinates
\[
(\xi, \eta) = (x, u(x,y)), \qquad (x,y) = (\xi, h(\xi,\eta)),
\]
which map every candidate fluid domain homeomorphically and smoothly onto the universal strip $S = \mathbb{R} \times (-1,1)$. Under this change of variables, the flow geometry is completely encoded by the height function $y = h(\xi,\eta)$, which represents the vertical position along the streamline $u = \eta$ at horizontal coordinate $\xi$. The harmonic equation $\Delta u = 0$ transforms into a  divergence-form elliptic equation for $h$:
\[
\operatorname{div}_{\xi,\eta} DF(\nabla_{\xi,\eta} h) = 0 \quad \text{in } S,
\qquad \text{where } F(z_1, z_2) = \frac{1+z_1^2}{2z_2}.
\]

Let $(u,\Gamma_1,\Gamma_2,\lambda,\theta)$ and $(\widetilde u,\widetilde \Gamma_1,\widetilde \Gamma_2,\widetilde{\lambda},\widetilde\theta)$ be two candidate jets issuing from the same nozzle, with downstream asymptotic slopes $\tan\theta, \tan\widetilde{\theta}$ and intercepts $b_i, \widetilde b_i$ ($i=1,2$). Denoting their height functions by $h$ and $\widetilde h$, their difference $w = h - \widetilde h$ satisfies a divergence-form elliptic equation
\[
\operatorname{div}(\mathcal{A} \nabla w) = 0 \quad \text{in } S,
\]
where the coefficient matrix $\mathcal{A}(Z)$ is symmetric and locally uniformly elliptic. On the upstream solid walls $\{\xi \le 0\}$, because the two flows share the identical nozzle geometry, $w$ satisfies the homogeneous Dirichlet conditions $w(\xi, \pm 1) = 0$. On the downstream free boundaries $\mathcal{T}_1 = \{(\xi,1)\mid \xi>0\}$ and $\mathcal{T}_2 = \{(\xi,-1)\mid \xi>0\}$, the constant-speed free-boundary conditions translate into conormal boundary conditions:
\[
\mathcal{A}\nabla w \cdot \nu_1 = -\frac{\lambda^2 - \widetilde\lambda^2}{2} \quad \text{on } \mathcal{T}_1,
\qquad
\mathcal{A}\nabla w \cdot \nu_2 = \frac{\lambda^2 - \widetilde\lambda^2}{2} \quad \text{on } \mathcal{T}_2,
\]
where $\nu_1 = (0, 1)$ and $\nu_2 = (0, -1)$ are the outward unit normals to $S$. Moreover, the upstream coincidence of the nozzle enforces the uniform asymptotic decay $w(\xi,\eta) \to 0$ as $\xi \to -\infty$.

The proof of complete uniqueness is carried out in two principal stages.

\paragraph{Stage 1.} We prove that the free-boundary speed is unique.

We first prove that $\lambda = \widetilde\lambda$. Arguing by contradiction, suppose without loss of generality that $\lambda > \widetilde\lambda$. Then the conormal boundary data on $\mathcal{T}_1$ and $\mathcal{T}_2$ are strictly negative and positive, respectively. Crucially, the $C^1$ regularity up to the detachment points and the smooth fit conditions enable us to compute the explicit normal derivatives at the detachment points:
\[
w_\eta(0,1) = \sqrt{1+g_1'(0)^2}\left(\frac{1}{\lambda} - \frac{1}{\widetilde\lambda}\right) < 0,
\qquad
w_\eta(0,-1) = \sqrt{1+g_2'(0)^2}\left(\frac{1}{\lambda} - \frac{1}{\widetilde\lambda}\right) < 0.
\]
Since $w(0,\pm 1) = 0$, this forces $w$ to be strictly positive immediately below the upper detachment point $\widehat A = (0,1)$ and strictly negative immediately above the lower detachment point $\widehat B = (0,-1)$.

Let $U_+$ and $U_-$ denote the connected components of $\{w > 0\}$ and $\{w < 0\}$ containing these respective local sign regions. Utilizing the sign of the conormal derivatives on $\mathcal{T}_1 \cup \mathcal{T}_2$ and the Hopf boundary point lemma, we prove that $w$ cannot attain a local maximum on $\mathcal{T}_1$ or a local minimum on $\mathcal{T}_2$. This boundary maximum principle yields two geometric separation properties:
\begin{enumerate}[label=(\roman*)]
	\item The closure $\overline{U_+}$ cannot intersect the lower free boundary $\mathcal{T}_2$ at any point where $w > 0$, and symmetrically $\overline{U_-}$ cannot reach $\mathcal{T}_1$ where $w < 0$.
	\item Both connected components $U_+$ and $U_-$, as well as their superlevel/sublevel sets, must be unbounded downstream, reaching $\xi \to +\infty$.
\end{enumerate}
On the other hand, the downstream far-field expansion of the height functions guarantees that
\begin{equation}
w(\xi,\eta) = A_0\xi + B_0 + D_0\eta + o(1) \quad \text{uniformly for } \eta \in [-1,1] \text{ as } \xi \to +\infty,\label{a-10}
\end{equation}
with $A_0 = \tan\theta- \tan{\widetilde\theta}$ and $D_0 = \frac{1}{\lambda\cos\theta} - \frac{1}{\widetilde\lambda\cos\widetilde{\theta}}$. Because both positive and negative regions extend downstream to infinity, the linear growth coefficient must vanish, enforcing $A_0 = 0$ (hence $\theta = \widetilde\theta$). Furthermore, the non-reaching property of $\overline{U_+}$ onto $\mathcal{T}_2 \cap \{w > 0\}$ combined with the affine asymptotic behaviors requires  $D_0 \ge 0$. However, substituting $\theta = \widetilde\theta$ and $\lambda > \widetilde\lambda$ into the formula for $D_0$ yields $D_0 < 0$, yielding a contradiction. Thus, we must have $\lambda = \widetilde\lambda$.

\paragraph{Stage 2.} We prove that the height function is unique and therefore show that the asymmetric jet flow is unique. 

Once $\lambda=\widetilde\lambda$ has been established,
the conormal boundary conditions become
\begin{equation}\label{a11}
	\mathcal A\nabla w\cdot\nu=0
	\quad\text{on }\mathcal T_1\cup\mathcal T_2.
\end{equation}
We first establish an extremum principle for $w$
(Lemma~\ref{lem-extrema}). More precisely, if $w$ attains a
nonpositive global minimum or a nonnegative global maximum
at a finite point of $\overline S$, then $w\equiv0$.

To explain the argument, suppose that $w\not\equiv0$ and
that $w$ attains a nonpositive global minimum at
$Z_*\in\overline S$. The strong maximum principle excludes
$Z_*\in S$. If $Z_*\in\mathcal T_1\cup\mathcal T_2$, the
Hopf boundary point lemma contradicts the homogeneous
conormal condition \eqref{a11}. Thus $Z_*$ lies on a fixed
wall or at one of the endpoints. Since $w=0$ there, the
minimum value is zero. The same maximum principles and
\eqref{a11} then give
$w>0$ in $S\cup\mathcal T_1\cup\mathcal T_2$.

In particular, $w\ge0$ on $\overline S$ and
$w(\widehat A)=0$. For sufficiently small $r>0$, set
$D_r=(-r,r)\times(1-r,1)$. The Hopf boundary point lemma
at $(-r,1)$, together with the positivity of $w$ on the
bottom and right sides of $D_r$, gives
\[
w(\xi,\eta)\ge c_0(1-\eta)
\qquad\text{on }\partial D_r\cap\{\eta<1\}
\]
for some $c_0>0$.

We next define
$h_\varepsilon(\xi,\eta)
=h\bigl(\xi,\eta-\varepsilon(1-\eta)\bigr)$
on $\overline{D_r}$, where $\varepsilon>0$ is sufficiently
small. The function $h_\varepsilon$ satisfies the same height
equation as $h$ in $D_r$ and agrees with $h$ on the top
side of $D_r$. Moreover, the local $C^1$ regularity of $h$
gives $0\le h-h_\varepsilon\le C\varepsilon(1-\eta)$ on
$\overline{D_r}$, for some $C>0$. Consequently,
$h_\varepsilon\ge\widetilde h$ on $\partial D_r$ when
$\varepsilon>0$ is sufficiently small.

The local $C^1$ regularity and the positivity of
$(h_\varepsilon)_\eta$ and $\widetilde h_\eta$ ensure that
the equation for $h_\varepsilon-\widetilde h$ is uniformly
elliptic in $D_r$. The maximum principle therefore yields
$h_\varepsilon\ge\widetilde h$ in $\overline{D_r}$. Since
$h_\varepsilon(\widehat A)=\widetilde h(\widehat A)$, while
$\lambda=\widetilde\lambda$ gives
$h_\eta(\widehat A)=\widetilde h_\eta(\widehat A)>0$, we
obtain
\[
0\le
\widetilde h_\eta(\widehat A)
-(1+\varepsilon)h_\eta(\widehat A)
=-\varepsilon h_\eta(\widehat A)<0,
\]
a contradiction.  

We next combine Lemma~\ref{lem-extrema} with
\eqref{a-10}.
If $A_0>0$, then $w\to+\infty$ uniformly as
$\xi\to+\infty$. Since $w\to0$ uniformly as
$\xi\to-\infty$ and $w=0$ on the fixed walls,
$w$ must attain a nonpositive global minimum at a finite
point of $\overline S$.
Lemma~\ref{lem-extrema} would then give $w\equiv0$,
contradicting the downstream limit.
Interchanging the two jets excludes $A_0<0$.
Thus $A_0=0$, and
\eqref{difference-asymptotic-coefficients} gives
$\theta=\widetilde\theta$.
Together with $\lambda=\widetilde\lambda$,
this also gives $D_0=0$.
The same argument excludes $B_0>0$ and $B_0<0$,
so $B_0=0$.
Consequently, $w\to0$ uniformly at both ends of $S$.
If $w\not\equiv0$, this convergence forces $w$ to
attain a positive global maximum or a negative global
minimum at a finite point of $\overline S$.
Lemma~\ref{lem-extrema} therefore yields
$h=\widetilde h$ on $\overline S$,
as asserted in Proposition~\ref{prop-fixed-speed-height}.

Finally, inverting the partial hodograph transformation
gives coincidence of the fluid domains, free boundaries,
and stream functions.
Hence the velocity fields $(u_y,-u_x)$ and
$(\widetilde u_y,-\widetilde u_x)$ agree.
Bernoulli's law, together with the common free-boundary
speed $\lambda$ and atmospheric pressure
$p_{\mathrm{atm}}$, then gives equality of the pressure
fields, completing the proof of Theorem~\ref{thm-main}.

\section{The partial hodograph transformation and the height function}\label{sec-height}
In this section, we first prove that the solution $u$ is
strictly increasing in the vertical direction throughout the
flow domain. This monotonicity allows us to introduce the
partial hodograph transformation and represent the streamlines
by a height function on a fixed strip. We then establish the
regularity and far-field asymptotic behavior of the height function,
which will be used to compare two jets in the subsequent sections.

The next lemma proves strict vertical monotonicity of $u$ in the flow domain. This property allows us to use the value of  $u$ as a coordinate and to define the height function.

\begin{lemma}\label{lem-monotonicity}
Let $(u,\Gamma_1,\Gamma_2,\lambda,\theta)$ be a solution to asymmetric incompressible jet flows in
the sense of Definition~\ref{def-jet}, and let $\Omega$ be its flow domain defined in \eqref{jet-domain}. Then
we have
$$-1<u(x,y)<1\qquad\text{and}\qquad u_y(x,y)>0\quad\text{for every
$(x,y)\in\Omega$}.$$
In particular, for every $x\in\R$, the function
$y\mapsto u(x,y)$ is strictly increasing on $[a(x),b(x)]$.
\end{lemma}
\begin{proof}
	By \eqref{nozzle-upstream} and \eqref{jet-downstream}, and the continuity of the boundary functions $a, b$ defined in \eqref{jet-ab},  there exists a constant $d_0 > 0$ such that
	\[
	b(x) - a(x) \ge d_0 \qquad \text{for all } x \in \mathbb{R}.
	\]
	
For $0<t<d_0$, define
\begin{equation}\label{rev-vertical-difference}
	\Omega_t=\{(x,y):a(x)<y<b(x)-t\},\qquad
	\varphi_t(x,y)=u(x,y+t)-u(x,y).
\end{equation}
 The properties of $u$ in \eqref{jet-equations} then ensure that
	\[
	\varphi_t \in C^2(\Omega_t) \cap C(\overline{\Omega_t}),
	\qquad
	\Delta \varphi_t = 0 \quad \text{in } \Omega_t,
	\qquad
	|\varphi_t| \le 2 \quad \text{on } \overline{\Omega_t}.
	\]
	Regarding the boundary values of $\varphi_t$, since $-1 \le u \le 1$ in $\overline\Omega$ and $u=1$ (resp.\ $u=-1$) on the upper (resp.\ lower) boundary of $\Omega$, we find
	\[
		\varphi_t(x, a(x))  = u(x, a(x)+t) - (-1)  \ge 0,\quad
		\varphi_t(x, b(x)-t) = 1 - u(x, b(x)-t) \ge 0
	\]
	for every $x \in \mathbb{R}$. Hence, $\varphi_t \ge 0$ on $\partial\Omega_t$.
	
	To propagate this nonnegativity into the interior of the unbounded domain $\Omega_t$, we apply the maximum principle for unbounded domains (Lemma~\ref{lem-unbounded} in Appendix~\ref{app-max}). Since $a$ is continuous with $a(0) = -1$, there exists $\rho_0 \in (0, 1/4)$ such that $a(x) > -3/2$ for $|x| \le \rho_0$. Consequently, the closed ball $\overline{B_{\rho_0}((0,-2))}$ lies outside the flow domain, namely,
	\[
	\overline{B_{\rho_0}((0,-2))} \subset \mathbb{R}^2 \setminus \Omega \subset \mathbb{R}^2 \setminus \Omega_t.
	\]
	This provides the geometric condition required by Lemma~\ref{lem-unbounded}. Applying this maximum principle to $-\varphi_t$, which is bounded and satisfies $-\varphi_t \le 0$ on $\partial\Omega_t$, yields
	\[
	\varphi_t(x,y) = u(x,y+t) - u(x,y) \ge 0 \qquad \text{for all } (x,y) \in \Omega_t \text{ and } 0 < t < d_0.
	\]
	
	Given any $(x,y) \in \Omega$, choosing $0 < t < \min\{d_0, b(x)-y\}$ and taking the limit as $t \downarrow 0$ leads to
	\[
	u_y(x,y) = \lim_{t \downarrow 0} \frac{u(x,y+t) - u(x,y)}{t} \ge 0.
	\]
	Because $u$ is harmonic, interior regularity guarantees that $u_y \in C^\infty(\Omega)$ and $\Delta u_y = 0$ in $\Omega$. In view of the boundary value of $u$, the strong maximum principle then yields that
	\[
	u_y(x,y) > 0 \qquad \text{in } \Omega.
	\]
	For each fixed $x \in \mathbb{R}$, this strict positivity along with $u \in C([a(x), b(x)])$ shows that $y \mapsto u(x,y)$ is strictly increasing on $[a(x), b(x)]$. Finally, the boundary values $u(x,a(x)) = -1$ and $u(x,b(x)) = 1$ immediately imply that $-1 < u(x,y) < 1$ for all $(x,y) \in \Omega$.
\end{proof}

Next, we use a partial hodograph transformation to express the
unknown flow domain in terms of a height function on a
fixed strip. For a corresponding formulation in the
water-wave setting, see \cite[Section~2.1]{KL23}. Precisely, by Lemma~\ref{lem-monotonicity}, we introduce the partial
hodograph transformation
$$
\begin{aligned}
T:\Omega&\to S=\R\times(-1,1)\\
(x,y)&\mapsto (\xi,\eta):=(x,u(x,y))
\end{aligned}
$$
Since $u\in C(\bar\Omega)$, one can extend
$T$ continuously to $\overline\Omega$ with values in $\overline S$.
The upper and lower free boundaries are mapped onto
$\mathcal{T}_1 =\{(\xi,1)\mid \xi>0\}$ and $\mathcal{T}_2 =\{(\xi,-1)\mid \xi>0\}$, respectively, while the nozzle
endpoints $A=(0,1)$ and $B=(0,-1)$ are mapped onto
$\widehat A=(0,1)$ and $\widehat B=(0,-1)$.
In the next proposition, we define the height function as the second component of $T^{-1}$ and prove its regularity in $S$ and up to the indicated parts of $\partial S$.

\begin{proposition}
	\label{prop-inverse}
	Let $(u,\Gamma_1,\Gamma_2,\lambda,\theta)$ be a solution to asymmetric incompressible jet flows
	in the sense of Definition~\ref{def-jet}.
	Then $T:\Omega\to S$ is a $C^\infty$ diffeomorphism.
		Its inverse has the form
		$\mathcal H(\xi,\eta)=(\xi,h(\xi,\eta))$, where the height
		function $h$ satisfies
	\[
	h\in C^\infty(S)\cap C(\overline S)
	\cap C^{1,\alpha}_{\mathrm{loc}}(S\cup\mathcal{T}_1\cup\mathcal{T}_2),
	\qquad
	h_\eta>0\quad\text{in }S\cup\mathcal{T}_1\cup\mathcal{T}_2.
	\]
	Here $\alpha\in(0,1)$ is the exponent in condition \textup{(J2)}.
	Moreover,  $
	h\in C^1\!\left(
	\overline S\cap
	\bigl(B_{r_*}(\widehat A)\cup B_{r_*}(\widehat B)\bigr)
	\right)$ for some $r_*>0$.
\end{proposition}

\begin{proof}
	This proof is split into two steps.
	
\textbf{Step 1}. We define $h$ and prove its continuity on $\overline S$.

By Lemma \ref{lem-monotonicity}, for each fixed $\xi \in \mathbb{R}$,  the map $y \mapsto u(\xi, y)$ is strictly increasing on $(a(\xi), b(\xi))$ and maps $(a(\xi), b(\xi))$ homeomorphically onto $(-1, 1)$. Consequently, the map
\[
T \colon \Omega \to S, \qquad T(x, y) = (x, u(x, y))
\]
is a bijection, and its inverse is given by
\[
\mathcal{H} \colon S \to \Omega, \qquad \mathcal{H}(\xi, \eta) = (\xi, h(\xi, \eta)),
\]
where $h(\xi, \cdot)$ denotes the inverse of $u(\xi, \cdot)$ on $(-1, 1)$ for each $\xi \in \mathbb{R}$.

Since $u$ is harmonic in $\Omega$, we have $u \in C^\infty(\Omega)$, which implies that $T \in C^\infty(\Omega, S)$. The Jacobian matrix of $T$ is
\[
DT(x,y) =
\begin{pmatrix}
	1 & 0 \\
	u_x(x,y) & u_y(x,y)
\end{pmatrix},
\]
which satisfies $\det DT(x,y) = u_y(x,y) > 0$ for all $(x, y) \in \Omega$. Since $T$ is a smooth bijection with a non-vanishing Jacobian determinant everywhere in $\Omega$, the inverse function theorem implies that its inverse $\mathcal{H}$ is also smooth. Therefore, $T$ is a $C^\infty$-diffeomorphism from $\Omega$ onto $S$.

We now extend $h$ to the closure $\overline S = \mathbb{R} \times [-1, 1]$. Since $u(\xi, a(\xi)) = -1$ and $u(\xi, b(\xi)) = 1$, the strictly increasing continuous function $y \mapsto u(\xi, y)$ extends to a homeomorphism from $[a(\xi), b(\xi)]$ onto $[-1, 1]$. We can therefore uniquely extend $h$ to $\overline S$ by setting
\[
h(\xi, -1) = a(\xi) \quad \text{and} \quad h(\xi, 1) = b(\xi) \qquad \text{for all } \xi \in \mathbb{R}.
\]
Under this definition, $h(\xi, \eta)$ is the unique root $y \in [a(\xi), b(\xi)]$ of the equation $u(\xi, y) = \eta$ for every $(\xi, \eta) \in \overline S$.

To verify that $h \in C(\overline S)$, consider a sequence $\{(\xi_n, \eta_n)\} \subset \overline S$ such that $(\xi_n, \eta_n) \to (\xi_0, \eta_0) \in \overline S$ as $n \to \infty$. By definition of $h$,
\[
a(\xi_n) \le h(\xi_n, \eta_n) \le b(\xi_n).
\]
The continuity of $a$ and $b$ implies that $\{h(\xi_n, \eta_n)\}_n$ is bounded. Let $\{h(\xi_{n_k}, \eta_{n_k})\}_k$ be any convergent subsequence with limit $y_*$. Passing to the limit in the inequalities yields
\[
a(\xi_0) \le y_* \le b(\xi_0).
\]
Moreover, by the continuity of $u$ on $\overline\Omega$,
\[
u(\xi_0, y_*) = \lim_{k \to \infty} u(\xi_{n_k}, h(\xi_{n_k}, \eta_{n_k})) = \lim_{k \to \infty} \eta_{n_k} = \eta_0.
\]
Since $y \mapsto u(\xi_0, y)$ is strictly increasing on $[a(\xi_0), b(\xi_0)]$, there exists a unique solution to $u(\xi_0, y) = \eta_0$ in this interval, which is precisely $h(\xi_0, \eta_0)$. Therefore, every convergent subsequence must converge to the same limit $y_* = h(\xi_0, \eta_0)$. It follows that the entire sequence satisfies $h(\xi_n, \eta_n) \to h(\xi_0, \eta_0)$ as $n \to \infty$, which proves that $h \in C(\overline S)$.
	
\textbf{Step 2}. We prove the regularity of $h$ up to the boundary.

We first consider the upper boundary $\{\eta=1\}$, while the lower boundary $\{\eta=-1\}$ is treated analogously. Differentiating the identity $u(x,f_1(x))=1$ with respect to $x$ gives
\[
u_x(x,f_1(x)) + f_1'(x) u_y(x,f_1(x)) = 0 \qquad \text{for every } x>0.
\]
Since $u_y > 0$ in $\Omega$ and $u_y$ is continuous up to $\Gamma_1$, we have $u_y \ge 0$ on $\Gamma_1$. Recalling that $\lvert\nabla u\rvert=\lambda$ on $\Gamma_1$, we obtain
\[
u_y(x,f_1(x)) = \frac{\lambda}{\sqrt{1+f_1'(x)^2}} > 0 \qquad \text{for every } x>0.
\]
Moreover, the continuity of $\nabla u$ at the detachment point $A$ in \eqref{jet-endpoint-regularity} and the smooth-fit condition $f_1'(0)=g_1'(0)$ in \eqref{jet-fit} ensure that this strictly positive lower bound extends continuously to $x=0$, so $u_y(A) > 0$.

Now let $P_0=(x_0,f_1(x_0))$ with $x_0>0$, and let
$\widehat P_0=(x_0,1)$. By condition \textup{(J2)}, $u$ is
of class $C^{1,\alpha}$ up to the $C^2$ boundary $\Gamma_1$
near $P_0$. Extend $u$ locally across $\Gamma_1$ as a
$C^{1,\alpha}$ function. Since $u_y(P_0)>0$, the implicit
function theorem gives a local $C^{1,\alpha}$ solution of
$u(\xi,y)=\eta$ near $\widehat P_0$. By the uniqueness proved
in Step~1, this solution is $y=h(\xi,\eta)$. Hence
\begin{equation*}
    h\in C^{1,\alpha}_{\mathrm{loc}}(S\cup\mathcal{T}_1),
    \qquad
    h_\eta=\frac{1}{u_y\circ\mathcal H}>0
    \quad\text{on }S\cup\mathcal{T}_1.
\end{equation*}

At $A$, the joined boundary $N_1\cup\{A\}\cup\Gamma_1$
is of class $C^1$ by \eqref{jet-fit}, and $u$ is of class
$C^1$ by \eqref{jet-endpoint-regularity}. Since $u_y(A)>0$,
the same extension and inverse-function argument, now in class
$C^1$, gives $h\in C^1(\overline S\cap B_{r_A}(\widehat A))$
for some $r_A>0$. The corresponding argument at the lower
boundary gives
$h\in C^{1,\alpha}_{\mathrm{loc}}(S\cup\mathcal{T}_2)$,
$h_\eta>0$ on $\mathcal{T}_2$, and
$h\in C^1(\overline S\cap B_{r_B}(\widehat B))$
for some $r_B>0$.

Setting $r_* = \min\{r_A, r_B\} > 0$, we conclude that $h \in C^1\!\left(\overline S \cap \bigl(B_{r_*}(\widehat A)\cup B_{r_*}(\widehat B)\bigr)\right)$, which completes the proof.
\end{proof}
\begin{remark}
	Differentiating $u(\xi,h(\xi,\eta))=\eta$ in $S$ and extending
	the resulting identities to the boundary by continuity, we obtain
	\[
	h_\eta=\frac{1}{u_y\circ\mathcal H}
	\quad\text{on }S\cup\mathcal{T}_1\cup\mathcal{T}_2
	\cup\{\widehat A,\widehat B\}.
	\]
	Since
	$u_y(x,f_i(x))=\lambda/\sqrt{1+f_i'(x)^2}$ for $x>0$,
	letting $x\downarrow0$ and using the continuity of $\nabla u$ (see \eqref{jet-endpoint-regularity})
	and the smooth-fit condition (see \eqref{jet-fit}) yields
	\begin{equation}\label{inverse-endpoint-derivatives}
		h_\eta(0,1)=\frac{\sqrt{1+g_1'(0)^2}}{\lambda},
		\qquad
		h_\eta(0,-1)=\frac{\sqrt{1+g_2'(0)^2}}{\lambda}.
	\end{equation}
\end{remark}

Finally, we establish the far-field asymptotic expansions for the stream
function $u$ and the height function $h$ at both the downstream and upstream ends.
These asymptotic formulas describe the convergence of the jet flow toward uniform parallel flows in the far fields and will be essential for the comparison of two candidate
jets as $\xi\to\pm\infty$.

\begin{lemma}\label{lem-farfield}
	Let $(u,\Gamma_1,\Gamma_2,\lambda,\theta)$ be a solution to asymmetric incompressible jet flows in the
	sense of Definition~\ref{def-jet}, and let $h$ be its height function.
	With $b_1,b_2$ as in \eqref{jet-downstream}, we have
	\begin{equation}\label{farfield-stream}
		\lim_{x\to+\infty}
		\sup_{a(x)\le y\le b(x)}
		\left|
		u(x,y)+1-\frac{2}{b_1-b_2}(y- x\tan\theta-b_2)
		\right|=0.
	\end{equation}
	Furthermore, at the downstream end, the height function satisfies
	\begin{equation}\label{farfield-height-downstream}
		\lim_{\xi\to+\infty}
		\sup_{-1\le\eta\le1}
		\left|
		h(\xi,\eta)-\xi\tan\theta-\frac{b_1+b_2}{2}
		-\frac{1}{\lambda\cos\theta}\eta
		\right|=0.
	\end{equation}
	At the upstream end, the height function satisfies
	\begin{equation*}
		\lim_{\xi\to-\infty}
		\sup_{-1\le\eta\le1}
		\left|
		h(\xi,\eta)-\frac{H_1+H_2}{2}
		-\frac{H_1-H_2}{2}\eta
		\right|=0.
	\end{equation*}
\end{lemma}
\begin{proof}
	The proof is divided into three steps.
	
	\textbf{Step 1}. We estimate the difference between $u$ and the
	affine function in \eqref{farfield-stream}.
	
	Let $W=b_1-b_2>0$ denote the asymptotic vertical width of the
	jet, and define
	\[
	z(x,y)
	=u(x,y)-\left[-1+\frac{2}{W}
	\bigl(y-x\tan\theta-b_2\bigr)\right].
	\]
	Fix $\varepsilon\in(0,W/4)$. By the downstream asymptotic
	behavior in \eqref{jet-downstream}, there exists $R>0$ such that
	\begin{equation}\label{farfield-boundary-errors}
		|a(x)-x\tan\theta-b_2|\le\varepsilon,
		\qquad
		|b(x)-x\tan\theta-b_1|\le\varepsilon
		\qquad\text{for all }x\ge R.
	\end{equation}
	For this choice of $R$, define
	\[
	D_R=\{(x,y)\in\Omega:x>R\}.
	\]
	We first estimate $z$ on $\partial D_R\cap\partial\Omega$.
	Since $u(x,a(x))=-1$, $u(x,b(x))=1$, and $b_1-b_2=W$,
	the inequalities \eqref{farfield-boundary-errors} give
$$
		|z(x,a(x))|
	 =\frac{2}{W}|a(x)-x\tan\theta-b_2|
		\le\frac{2\varepsilon}{W},\quad
		|z(x,b(x))|
		=\frac{2}{W}|b(x)-x\tan\theta-b_1|
		\le\frac{2\varepsilon}{W}
	$$
	for every $x\ge R$.
	
	Next, we estimate $z$ throughout $\overline{D_R}$.
	For every $(x,y)\in\overline{D_R}$, we have
	$a(x)\le y\le b(x)$. Together with
	\eqref{farfield-boundary-errors}, this implies
	\[
	-\varepsilon
	\le y-x\tan\theta-b_2
	\le W+\varepsilon
	\qquad\text{in }\overline{D_R}.
	\]
	Using the definition of $z$ and the bound
	$-1\le u\le1$ in $\overline\Omega$, we therefore obtain
	\[
	|z(x,y)|
	\le 2+\frac{2\varepsilon}{W}
	\le\frac52
	\qquad\text{for every }(x,y)\in\overline{D_R}.
	\]
	Finally, the definition of $z$ and the equation
	$\Delta u=0$ in $\Omega$ show that
	$z\in C(\overline{D_R})\cap C^2(D_R)$ satisfies
	\begin{equation}\label{b8}
		\Delta z=0
		\quad\text{in }D_R,
		\qquad
		|z|\le\frac{2\varepsilon}{W}
		\quad\text{on }\partial D_R\cap\partial\Omega.
	\end{equation}

	\textbf{Step 2}. We bound $|z|$ by applying the maximum principle.
	
	To control $|z|$ on the unbounded domain $D_R$, we construct a positive barrier $\Phi_R$ satisfying
	\begin{equation*}
		\Delta \Phi_R = 0 \quad \text{in } D_R, \qquad
		|z| \le \frac{2\varepsilon}{W} + \Phi_R \quad \text{on } \partial D_R, \qquad
		\lim_{x \to +\infty} \sup_{a(x) \le y \le b(x)} \Phi_R(x, y) = 0.
	\end{equation*}
	Let $c_0 = (b_1 + b_2)/2$. We consider the function
	\[
	\Phi_R(x, y) = C e^{-\kappa[(x - R) + \tan\theta(y - R\tan\theta- c_0)]} \cos\!\bigl(\kappa(y -  x\tan\theta - c_0)\bigr),
	\]
	where $\kappa= 1/W$ and the constant $C > 0$ will be determined later. Direct differentiation shows that
	\[
	\Delta\Phi_R=\kappa^2\bigl[(1+\tan^2\theta)-(1+\tan^2\theta)\bigr]\Phi_R=0\quad\text{in }D_R.
	\]

	By \eqref{farfield-boundary-errors}, for any $(x, y) \in \overline{D_R}$ we have
	\begin{equation}\label{b-10}
		|y -  x\tan\theta - c_0| \le \frac{W}{2} + \varepsilon < \frac{3W}{4},
	\end{equation}
	which ensures that
	\[
	\frac{1}{2} < \cos(3/4) \le \cos\!\left(\frac{y -  x\tan\theta - c_0}{W}\right) \le 1 \qquad \text{in}\ \  \overline{D_R}.
	\]
	Hence, for $x = R$ and $a(R) \le y \le b(R)$, using \eqref{b-10}, we obtain
	$$
	\begin{aligned}
	\Phi_R(R, y)& = C e^{-\kappa\tan\theta(y - R\tan\theta  - c_0)} \cos\kappa(y - R\tan\theta - c_0)\\
	&\ge \frac{C}{2} e^{-\kappa|\tan\theta|\,|y - R\tan\theta - c_0|}\\
	&\ge \frac{C}{2} e^{-3|\tan\theta|/4}.
	\end{aligned}
    $$
	Since $|z(R, y)| \le 5/2$, choosing $C = 5 e^{3|\tan\theta|/4}$ ensures that
	$$\Phi_R(R, y) \ge 5/2 \ge |z(R, y)|\quad\text{on}~\{x = R\} \cap \overline\Omega.$$ Combining this with the boundary estimate in \eqref{b8}, we deduce that
	$$|z| \le \frac{2\varepsilon}{W} + \Phi_R\quad\text{on}~\partial D_R.$$
	Furthermore, for any $(x, y) \in \overline{D_R}$, using \eqref{b-10} again, we have
	\[
	\begin{aligned}
		0 < \Phi_R(x, y) &= 5e^{3|\tan\theta|/4} e^{-\kappa\sec^2\theta(x - R)} e^{-\kappa\tan\theta(y -  x\tan\theta - c_0)} \cos(\kappa(y -  x\tan\theta - c_0)) \\
		&\le 5 e^{3|\tan\theta|/2} e^{-\kappa\sec^2\theta(x - R)}\to0,\qquad\text{as}~x\to+\infty.
	\end{aligned}
	\]
	 Since the complement $\mathbb{R}^2 \setminus D_R$ contains a ball (for instance $\overline{B_1((R-2,0))}$), applying the maximum principle on unbounded domains (Lemma~\ref{lem-unbounded}) to $\pm z - \frac{2\varepsilon}{W} - \Phi_R$ yields
	\[
	|z(x, y)| \le \frac{2\varepsilon}{W} + 5e^{3|\tan\theta|/2} e^{-(1 + \tan^2\theta)(x - R)/W} \qquad \text{in } \overline{D_R}.
	\]
	Taking the supremum over $y \in [a(x), b(x)]$ and sending $x \to +\infty$, we find
	\[
	\limsup_{x \to +\infty} \sup_{a(x) \le y \le b(x)} |z(x, y)| \le \frac{2\varepsilon}{W}.
	\]
	Letting $\varepsilon \downarrow 0$ establishes the asymptotic formula \eqref{farfield-stream}.
	
	\textbf{Step 3}. We derive the asymptotic formulas for the height function $h$.
	
	By Lemma~\ref{lem-monotonicity}, for each fixed $\xi \in \mathbb{R}$, the mapping $y \mapsto u(\xi, y)$ is a strictly increasing homeomorphism from $[a(\xi), b(\xi)]$ onto $[-1, 1]$. The height function defined in Proposition~\ref{prop-inverse} satisfies the identity
	\[
	u(\xi, h(\xi, \eta)) = \eta, \qquad a(\xi) \le h(\xi, \eta) \le b(\xi) \qquad \text{for all } (\xi, \eta) \in \mathbb{R} \times [-1, 1].
	\]
	Substituting $y = h(\xi, \eta)$ into the definition of $z$, we obtain
	\[
	z(\xi, h(\xi, \eta)) = \eta - \left[ -1 + \frac{2}{W}\bigl(h(\xi, \eta) - \xi\tan\theta - b_2\bigr) \right]
	= \eta + 1 - \frac{2}{W}\bigl(h(\xi, \eta) - \xi\tan\theta - b_2\bigr).
	\]
	 Using $b_2 + \frac{W}{2} = \frac{b_1 + b_2}{2}$, it follows that
	\[
	h(\xi, \eta) - \xi\tan\theta - \frac{b_1 + b_2}{2} - \frac{W}{2}\eta = -\frac{W}{2} z(\xi, h(\xi, \eta)).
	\]
	Consequently, using \eqref{farfield-stream}, we have
	\[
	\sup_{-1 \le \eta \le 1} \left| h(\xi, \eta) - \xi\tan\theta - \frac{b_1 + b_2}{2} - \frac{W}{2}\eta \right| \le \frac{W}{2} \sup_{a(\xi) \le y \le b(\xi)} |z(\xi, y)| \to 0 \qquad \text{as } \xi \to +\infty.
	\]
	From the downstream flux relation \eqref{jet-flux-width}, we have $\frac{W}{2} = \frac{b_1 - b_2}{2} = \frac{1}{\lambda\cos\theta}$. Substituting this into the above limit proves \eqref{farfield-height-downstream}.
	
	Finally, one can apply the same argument to obtain the upstream asymptotic formula for $h$. This completes the proof.
\end{proof}

\section{The elliptic equations of the height functions}
\label{sec-equation}

In this section, we formulate the framework for comparing two candidate jet flows on the reference strip. We first derive the nonlinear divergence-form elliptic equation satisfied by the height function of an individual jet flow. Next, considering two candidate asymmetric jets issuing from the same nozzle with identical mass flux, we linearize the problem and establish the elliptic equation, the conormal boundary conditions, and the far-field asymptotics satisfied by the difference of their height functions. Finally, we verify the local uniform ellipticity of the linearized operator, which serves as the basis for the comparison and maximum-principle arguments in Sections~\ref{sec-speed} and \ref{sec-main-proof}.

The nonlinear equation governing the height function follows directly from the harmonicity of the physical stream function through the partial hodograph transformation.

\begin{proposition}
	\label{prop-height-equation}
	Let $(u,\Gamma_1,\Gamma_2,\lambda,\theta)$ be a solution to asymmetric incompressible jet flows in the sense of Definition~\ref{def-jet}, and let $h$ be its associated height function established in Proposition~\ref{prop-inverse}. Then $h$ satisfies
	\begin{equation}\label{equation-height}
		\partial_\xi\left(\frac{h_\xi}{h_\eta}\right)
		-\partial_\eta\left(\frac{1+h_\xi^2}{2h_\eta^2}\right)=0
		\qquad\text{in }S.
	\end{equation}
\end{proposition}

\begin{proof}
	Let $X=(x,y)\in\Omega$, $Z=(\xi,\eta)\in S$, and let $\mathcal H(Z)=(\xi,h(\xi,\eta))$ be the inverse transformation. By definition of the height function,
	\[
	u(\xi,h(\xi,\eta))=\eta\qquad\text{for all }(\xi,\eta)\in S.
	\]
	Differentiating this identity with respect to $\xi$ and $\eta$ yields
	\[
	(u_x\circ\mathcal H)+(u_y\circ\mathcal H)h_\xi=0,
	\qquad
	(u_y\circ\mathcal H)h_\eta=1.
	\]
	Since $h_\eta>0$ in $S$ by Proposition~\ref{prop-inverse}, we can express the physical gradient $(\nabla_Xu)\circ\mathcal H$ in terms of the derivatives of $h$:
	\begin{equation}\label{equation-gradient-matrix}
		(\nabla_Xu)\circ\mathcal H=
		\begin{pmatrix}
			-h_\xi/h_\eta\\[1mm]
			1/h_\eta
		\end{pmatrix}
		\qquad\text{in }S.
	\end{equation}
	
	Next, we calculate the second derivatives of $u$ in hodograph coordinates. Applying the chain rule to the components of $\nabla_X u$ yields
	\begin{align*}
		\partial_\xi(u_x\circ\mathcal H)
		&=(u_{xx}\circ\mathcal H)+h_\xi(u_{xy}\circ\mathcal H),
		\qquad
		\partial_\eta(u_x\circ\mathcal H)
		=h_\eta(u_{xy}\circ\mathcal H),\\[1mm]
		\partial_\eta(u_y\circ\mathcal H)
		&=h_\eta(u_{yy}\circ\mathcal H).
	\end{align*}
	Solving for the spatial second derivatives gives
	\[
	u_{xx}\circ\mathcal H
	=\partial_\xi(u_x\circ\mathcal H)
	-\frac{h_\xi}{h_\eta}\partial_\eta(u_x\circ\mathcal H),
	\qquad
	u_{yy}\circ\mathcal H
	=\frac{1}{h_\eta}\partial_\eta(u_y\circ\mathcal H).
	\]
	Substituting \eqref{equation-gradient-matrix} into these relations and using $\Delta u = u_{xx}+u_{yy}=0$ in $\Omega$, we obtain
	\begin{align*}
		0=(\Delta u)\circ\mathcal H=-\partial_\xi\left(\frac{h_\xi}{h_\eta}\right)
		+\partial_\eta\left(\frac{1+h_\xi^2}{2h_\eta^2}\right).
	\end{align*}
	This gives \eqref{equation-height}.
\end{proof}

To prepare for the comparison of two candidate height functions, we express \eqref{equation-height} in divergence form. Define the function $F:\mathbb R\times(0,\infty)\to\mathbb R$ by $
F(z_1,z_2)=\frac{1+z_1^2}{2z_2}$.
Then a direct computation gives
\begin{equation}
DF(\nabla h)=
\begin{pmatrix}
	h_\xi/h_\eta\\[1.5mm]
	-(1+h_\xi^2)/(2h_\eta^2)
\end{pmatrix}.\label{equation-DF}
\end{equation}
Consequently, equation \eqref{equation-height} is equivalent to the divergence-form equation
\begin{equation}\label{equation-height-F}
	\operatorname{div}_{\xi,\eta}DF(\nabla_{\xi,\eta}h)=0
	\qquad\text{in }S.
\end{equation}

Let $(u,\Gamma_1,\Gamma_2,\lambda,\theta)$ and $(\tilde u,\tilde\Gamma_1,\tilde\Gamma_2,\tilde\lambda,\tilde\theta)$  be two solutions to asymmetric incompressible jets issuing from the same prescribed nozzle. Let $h,\tilde h:\overline S\to\mathbb R$ be their respective height functions.
By decreasing the two H\"older exponents if necessary, we use a common
exponent $\alpha\in(0,1)$ for both solutions.
We define the difference function $
w=h-\tilde h$,
and consider the family
\[
h_t=(1-t)\tilde h+th
\qquad\text{for }t\in[0,1].
\]
Since $(h_t)_\eta=(1-t)\tilde h_\eta+th_\eta>0$ in $S\cup\mathcal{T}_1\cup\mathcal{T}_2$ for each $t\in[0,1]$, the linearized coefficient matrix given by
\begin{equation}\label{difference-matrix}
	\mathcal A(Z)= \int_0^1 D^2F(\nabla h_t(Z))\,dt \qquad\text{for }Z\in S\cup\mathcal{T}_1\cup\mathcal{T}_2
\end{equation}
is well-defined.

Let $(\tan\theta,b_1,b_2)$ and $(\tan\tilde\theta,\tilde b_1,\tilde b_2)$ denote the downstream asymptotic parameters of the two flows. We introduce the asymptotic difference coefficients
\begin{equation}\label{difference-asymptotic-coefficients}
		A_0=\tan\theta-\tan\tilde\theta,\quad
		B_0=\frac{b_1+b_2}{2}-\frac{\tilde b_1+\tilde b_2}{2},\quad
		D_0=\frac{\sec\theta}{\lambda}-\frac{\sec\tilde\theta}{\tilde\lambda}.
\end{equation}

The following proposition establishes the mixed boundary value problem and asymptotic properties satisfied by $w$ on the fixed strip $S$.

\begin{proposition}\label{prop-height-difference}
	The height-function difference $
	w\in C(\overline S)\cap C^\infty(S)
    \cap C^{1,\alpha}_{\mathrm{loc}}(S\cup\mathcal{T}_1\cup\mathcal{T}_2)$
	satisfies the  boundary value problem
	\begin{equation}\label{difference-problem}
        \left\{
        \begin{aligned}
            &\operatorname{div}(\mathcal A\nabla w)=0
            &&\text{in }S,\\[1mm]
            &w(\xi,-1)=w(\xi,1)=0
            &&\text{for }\xi\le0,\\[1mm]
            &\mathcal A\nabla w\cdot\nu_1
            =-\frac{\lambda^2-\tilde\lambda^2}{2} \ \text{and}\   w=f_1(\xi)-\tilde{f}_1(\xi)
            &&\text{on }\mathcal{T}_1,\\[1mm]
            &\mathcal A\nabla w\cdot\nu_2
            =\frac{\lambda^2-\tilde\lambda^2}{2}  \ \text{and} \ w=f_2(\xi)-\tilde{f}_2(\xi)
            &&\text{on }\mathcal{T}_2.
        \end{aligned}
        \right.
    \end{equation}
	where $\nu_1 = (0,1)$ on $\mathcal{T}_1$ and $\nu_2 = (0,-1)$ on $\mathcal{T}_2$ denote the outward unit normal vectors. Furthermore, $w$ fulfills the far-field asymptotic behavior
	\begin{equation}\label{difference-asymptotics}
			 \lim_{\xi\to-\infty}
			\sup_{-1\le\eta\le1}|w(\xi,\eta)|=0,\qquad
	\lim_{\xi\to+\infty}
			\sup_{-1\le\eta\le1}
			\left|w(\xi,\eta)-A_0\xi-B_0-D_0\eta\right|=0.
	\end{equation}
	The matrix $\mathcal A$ is symmetric and satisfies
    \begin{equation}\label{difference-coefficient-regularity}
        \mathcal A\in C^\infty(S;\mathbb R^{2\times2})
        \cap C^{0,\alpha}_{\mathrm{loc}}
        (S\cup\mathcal{T}_1\cup\mathcal{T}_2;\mathbb R^{2\times2}).
    \end{equation}
    In particular, its coefficients are H\"older continuous up to
    the transformed free boundaries away from $\widehat A$ and
    $\widehat B$. For every compact set $K\subset S\cup\mathcal{T}_1\cup\mathcal{T}_2$, there exist positive constants $\mu_K$ and $\Lambda_K$ such that
	\begin{equation}\label{difference-ellipticity}
		\mu_K|\mathbf v|^2
		\le\mathbf v^T\mathcal A(Z)\mathbf v
		\le\Lambda_K|\mathbf v|^2
		\qquad\text{for all }Z\in K\text{ and }\mathbf v\in\mathbb R^2.
	\end{equation}
\end{proposition}

\begin{proof}
	We divide the proof into four steps.
	
	\textbf{Step 1}. We derive the equation for $w$.
	
	By Proposition~\ref{prop-inverse}, $h,\tilde h\in C(\overline S)\cap C^\infty(S)\cap C^{1,\alpha}_{\mathrm{loc}}(S\cup\mathcal{T}_1\cup\mathcal{T}_2)$. Equation \eqref{equation-height-F} ensures that
	\[
	\operatorname{div}DF(\nabla h)=0,\qquad
	\operatorname{div}DF(\nabla\tilde h)=0
	\qquad\text{in }S.
	\]
Thus we have
	\[
DF(\nabla h)-DF(\nabla\tilde h)
		=\left(\int_0^1 D^2F(\nabla h_t)\,dt\right)(\nabla h-\nabla\tilde h)
		=\mathcal A\nabla w.
	\]
	Taking the divergence on both sides yields $\operatorname{div}(\mathcal A\nabla w)=0$ in $S$.
	
	\textbf{Step 2}. We verify the boundary conditions.
	
	Since both flows share the same solid nozzle walls, we have
	\[
	h(\xi,-1)=\tilde h(\xi,-1)=g_2(\xi),\qquad
	h(\xi,1)=\tilde h(\xi,1)=g_1(\xi)
	\qquad\text{for }\xi\le0.
	\]
	Hence, $w(\xi,-1)=w(\xi,1)=0$ for all $\xi\le0$.
	
	Along the free boundaries $\mathcal{T}_1\cup\mathcal{T}_2$, formula \eqref{equation-gradient-matrix} applies up to the boundary for both candidate flows, which implies
	\[
	\frac{1+h_\xi^2}{h_\eta^2}=\lambda^2,\qquad
	\frac{1+\tilde h_\xi^2}{\tilde h_\eta^2}=\tilde\lambda^2
	\qquad\text{on }\mathcal{T}_1\cup\mathcal{T}_2.
	\]
	Using \eqref{equation-DF} and the definition of $\mathcal{A}\nabla w$, we have
	\[
	(\mathcal A\nabla w)_2
	=-\frac{1+h_\xi^2}{2h_\eta^2}+\frac{1+\tilde h_\xi^2}{2\tilde h_\eta^2}
	=-\frac{\lambda^2-\tilde\lambda^2}{2}
	\qquad\text{on }\mathcal{T}_1\cup\mathcal{T}_2.
	\]
	Recalling that the outward unit normals are $\nu_1=(0,1)$ on $\mathcal{T}_1$ and $\nu_2=(0,-1)$ on $\mathcal{T}_2$,  \eqref{difference-problem} follows.
	
	\textbf{Step 3}. We prove the far-field asymptotic expansions.
	
	Applying Lemma~\ref{lem-farfield} to both $h$ and $\tilde h$ yields the upstream limit
	\[
	\lim_{\xi\to-\infty}
	\sup_{-1\le\eta\le1}
	\left|h(\xi,\eta)-\frac{H_1+H_2}{2}-\frac{H_1-H_2}{2}\eta\right|=0,
	\]
	together with the identical relation for $\tilde h$. The triangle inequality immediately gives the upstream decay in \eqref{difference-asymptotics}.
	
	Similarly, using Lemma~\ref{lem-farfield} and \eqref{difference-asymptotic-coefficients}, one can establish the downstream limit in \eqref{difference-asymptotics}.
	
\textbf{Step 4}. We prove the regularity and ellipticity of $\mathcal A$.

Recall from \eqref{difference-matrix} that the linearized coefficient matrix is given by
\[
\mathcal{A}(Z) = \int_0^1 D^2F(\nabla h_t(Z))\,dt, \qquad Z \in S \cup \mathcal{T}_1 \cup \mathcal{T}_2.
\]
By the definition of $F$, a direct computation gives that the Hessian matrix of $F$ is
\begin{equation}\label{difference-Hessian}
	D^2F(z_1,z_2) =
	\begin{pmatrix}
		z_2^{-1} & -z_1 z_2^{-2} \\[1.5mm]
		-z_1 z_2^{-2} & (1+z_1^2)z_2^{-3}
	\end{pmatrix}.
\end{equation}
Since $D^2F(z_1,z_2)$ is symmetric for $z_2>0$,
$\mathcal A(Z)$ is symmetric. The interior smoothness of
$h,\widetilde h$ and the positivity of $(h_t)_\eta$ imply
$\mathcal A\in C^\infty(S;\mathbb R^{2\times2})$.
To prove H\"older continuity up to the free boundaries, fix a
compact set $K\subset S\cup\mathcal{T}_1\cup\mathcal{T}_2$.
Proposition~\ref{prop-inverse} gives a positive lower bound for
$h_\eta,\widetilde h_\eta$ and finite $C^{1,\alpha}$ norms on
$K$. Since $D^2F$ is smooth for $z_2>0$, formula
\eqref{difference-matrix} gives a constant $C_K>0$ such that
\begin{equation*}
    \begin{aligned}
        \lVert\mathcal A(Z)-\mathcal A(Z')\rVert
        &\le C_K\bigl(|\nabla h(Z)-\nabla h(Z')|
        +|\nabla\widetilde h(Z)-\nabla\widetilde h(Z')|\bigr)\\
        &\le C_K|Z-Z'|^\alpha
        \qquad\text{for }Z,Z'\in K.
    \end{aligned}
\end{equation*}
This proves \eqref{difference-coefficient-regularity}.

Next, we establish the ellipticity of $\mathcal{A}$. For every $\mathbf{v} = (v_1, v_2)^T \in \mathbb{R}^2 \setminus \{\mathbf{0}\}$ and $z_2 > 0$, formula \eqref{difference-Hessian} gives
\[
\mathbf{v}^T D^2F(z_1,z_2)\mathbf{v} = \frac{(z_2 v_1 - z_1 v_2)^2 + v_2^2}{z_2^3} > 0,
\]
showing that $D^2F$ is strictly positive definite on $\mathbb{R} \times (0,\infty)$.

Given any compact set $K \subset S \cup \mathcal{T}_1 \cup \mathcal{T}_2$, the continuity of $\nabla h, \nabla \tilde{h}$ and the condition $h_\eta, \tilde{h}_\eta > 0$ imply that $(h_t)_\eta \ge c_K > 0$ on $K \times [0,1]$. Hence, $\mathcal{K} := \{\nabla h_t(Z) : Z \in K,\, t \in [0,1]\}$ is a compact subset of $\mathbb{R} \times (0,\infty)$. By the strict positivity and continuity of $D^2F$ on $\mathcal{K}$, there exist constants $0 < \mu_K \le \Lambda_K$ such that
\[
\mu_K |\mathbf{v}|^2 \le \mathbf{v}^T D^2F(\nabla h_t(Z))\mathbf{v} \le \Lambda_K |\mathbf{v}|^2 \qquad \text{for all } Z \in K,\, t \in [0,1],\, \mathbf{v} \in \mathbb{R}^2.
\]
Integrating over $t \in [0,1]$ yields
\[
\mu_K |\mathbf{v}|^2 \le \mathbf{v}^T \mathcal{A}(Z)\mathbf{v} = \int_0^1 \mathbf{v}^T D^2F(\nabla h_t(Z))\mathbf{v}\,dt \le \Lambda_K |\mathbf{v}|^2,
\]
which proves the local uniform ellipticity \eqref{difference-ellipticity}
and completes the proof.
\end{proof}

\section{The uniqueness of the speed \texorpdfstring{$\lambda$}{lambda} on the free boundaries}
\label{sec-speed}

In this section, we complete the first stage of the uniqueness proof by showing that any two candidate jets must have the same free-boundary speed $\lambda$.

\begin{theorem}\label{prop-speed}
	Let $(u,\Gamma_1,\Gamma_2,\lambda,\theta)$ and $(\tilde u,\tilde \Gamma_1,\tilde \Gamma_2,\tilde\lambda,\tilde\theta)$ be two solutions to asymmetric incompressible jet flows in the sense of Definition~\ref{def-jet}. Then their free-boundary speeds must coincide, namely, $
	\lambda=\tilde\lambda$.
\end{theorem}

Throughout this section, we let $h$ and $\tilde h$ denote the height functions of the two candidate jets as established in Proposition~\ref{prop-inverse}, and we set $w = h - \tilde h$. We argue by contradiction. If the two speeds were distinct, we may assume without loss of generality   that
\begin{equation}\label{standing-speed-assumption}
	\lambda > \tilde\lambda.
\end{equation}
  Under \eqref{standing-speed-assumption}, the conormal boundary conditions in \eqref{difference-problem} immediately yield
\begin{equation}\label{conormal-signs}
	\mathcal A\nabla w\cdot\nu_1 = -\frac{\lambda^2-\tilde\lambda^2}{2} < 0 \quad\text{on }\mathcal{T}_1,
	\qquad
	\mathcal A\nabla w\cdot\nu_2 = \frac{\lambda^2-\tilde\lambda^2}{2} > 0 \quad\text{on }\mathcal{T}_2.
\end{equation}

In subsection \ref{subsec-boundary-sign-regions}, our analysis will be carried out under the assumption that there exists a constant $\delta \in (0, 2/3)$ such that
\begin{equation}\label{comparison-endpoint-signs}
	w(0,\eta) > 0 \quad\text{for } 1-\delta < \eta < 1,
	\qquad
	w(0,\eta) < 0 \quad\text{for } -1 < \eta < -1+\delta.
\end{equation}
In fact, using \eqref{conormal-signs}, one can easily verify \eqref{comparison-endpoint-signs}, which will be rigorously established in subsection~\ref{s4.2} by exploiting the smooth fit conditions \eqref{jet-fit} at the endpoints $A$ and $B$ of the asymmetric nozzle walls.

Under the standing assumption \eqref{comparison-endpoint-signs}, we define the following two sets:
\begin{itemize}
	\item $U_+$: the connected component of the open set $\{Z\in S : w(Z) > 0\}$ containing the vertical line segment $\{0\}\times(1-\delta,1)$;
	\item $U_-$: the connected component of the open set $\{Z\in S : w(Z) < 0\}$ containing the vertical line segment $\{0\}\times(-1,-1+\delta)$.
\end{itemize}

\subsection{Properties of the sets $U_+$ and $U_-$}
\label{subsec-boundary-sign-regions}

Throughout this subsection, we assume that \eqref{standing-speed-assumption} and \eqref{comparison-endpoint-signs} hold. We first show that the definite signs of the conormal derivatives in \eqref{conormal-signs} prevent $w$ from attaining a local maximum on $\mathcal{T}_1$ or a local minimum on $\mathcal{T}_2$. Using this boundary maximum principle, we then prove that the closure $\overline{U_+}$ cannot intersect the lower free boundary $\mathcal{T}_2$ where $w > 0$, nor can $\overline{U_-}$ intersect the upper free boundary $\mathcal{T}_1$ where $w < 0$.

The following lemma establishes a boundary maximum principle that will be applied throughout the downstream analysis.

\begin{lemma}\label{strip-boundary-extrema}
	The difference function $w$ cannot attain a local maximum relative to $\overline S$ at any point of $\mathcal{T}_1$, nor can it attain a local minimum relative to $\overline S$ at any point of $\mathcal{T}_2$.
\end{lemma}
\begin{proof}
	Suppose by contradiction that $Z\in \mathcal{T}_1$ is a local maximum point of $w$ relative to $\overline S$.
	Since $w$ is $C^1$ up to $\mathcal{T}_1$, we have
	\begin{equation}\label{rev-boundary-maximum-derivatives}
		w_\xi(Z)=0,\qquad
		w_\eta(Z)=\lim_{t\downarrow0}
		\frac{w(Z)-w(Z-t\nu_1)}{t}\ge0.
	\end{equation}
	Using \eqref{rev-boundary-maximum-derivatives} and $\mathcal A_{22}(Z)>0$, we obtain
	\begin{equation}\label{rev-boundary-maximum-conormal}
		\mathcal A(Z)\nabla w(Z)\cdot\nu_1
		=\mathcal A_{22}(Z)w_\eta(Z)\ge0,
	\end{equation}
	contrary to \eqref{conormal-signs}, where $\nu_1=(0,1)$ is the unit normal vector at $Z$.
	The impossibility of a local minimum on $\mathcal{T}_2$ follows by an entirely symmetric argument applied to $-w$, and we omit the details.
\end{proof}

The next lemma shows that the positive component $U_+$ initiating near the upper detachment point $\widehat A$ remains separated from the lower free boundary $\mathcal{T}_2$ on $\{w > 0\}$.

\begin{lemma}\label{strip-components}
	The connected components $U_+$ and $U_-$ satisfy the disjointness relations:
	\[
	\overline{U_+}\cap\{Z\in\mathcal{T}_2:w(Z)>0\}=\varnothing,
	\qquad
	\overline{U_-}\cap\{Z\in\mathcal{T}_1:w(Z)<0\}=\varnothing.
	\]
\end{lemma}
\begin{proof}
	It suffices to prove the assertion for $U_+$, as the corresponding property for $U_-$ is established analogously. Suppose, by contradiction, that there exists a point $Z_r=(r,-1) \in \overline{U_+}$ for some $r > 0$ such that $w(Z_r) > 0$. By the continuity of $w$, we can choose $\rho_0 \in (0, r)$ such that
	\[
	B_{\rho_0}(Z_r)\cap\overline S \subset \{w > 0\},
	\qquad
	B_{\rho_0}(Z_r)\cap S \subset U_+.
	\]
	Fix $P_+ = (0,1-\delta/2) \in U_+$ and $P_r = (r,-1+\eta_r) \in U_+$, where $0 < \eta_r < \min\{\rho_0, 1\}$. Because $U_+$ is a connected open domain in $\mathbb R^2$, it is path-connected. Hence, there exists a simple polygonal path in $U_+$ joining $P_+$ to $P_r$. By connecting this path with the vertical segment from $\widehat A=(0,1)$ to $P_+$ and the vertical segment from $P_r$ to $Z_r$, and removing the loops, we can obtain an injective piecewise linear path $\gamma \colon [0,1] \to \overline S$ satisfying
	\begin{equation}\label{comparison-first-arc}
		\gamma(0)=\widehat A, \qquad
		\gamma(1)=(r,-1), \qquad
		\gamma((0,1))\subset U_+, \qquad
		w\ge0\quad\text{on }\gamma([0,1]).
	\end{equation}
	
	Since $\gamma$ is compact, there exists $L_0 > 0$ such that $\gamma \subset \{(\xi,\eta) \in \overline S : \xi > -L_0\}$. For each $L \ge L_0$, consider the closed polygonal curve
	\[
	J_L = \gamma \cup \bigl([-L,0]\times\{1\}\bigr) \cup \bigl(\{-L\}\times[-1,1]\bigr) \cup \bigl([-L,r]\times\{-1\}\bigr).
	\]
	Since $\gamma$ meets $\partial S$ only at $\widehat A$ and $(r,-1)$, $J_L$ is a simple closed curve enclosing a bounded domain $G_L \subset S$ with $\overline{G_L} \subset \overline S$ (see Figure~\ref{fig-first-comparison-domain}).
	
	\begin{figure}[htbp]
		\centering
		\includegraphics[width=50mm]{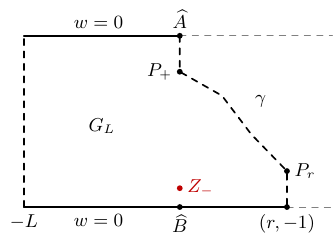}
		\caption{The path $\gamma$ and the domain $G_L$.}
		\label{fig-first-comparison-domain}
	\end{figure}
	
	For any $L \ge L_0$,  the upstream asymptotic decay \eqref{difference-asymptotics} of $w$ implies that
	\begin{equation*}
		\varepsilon_L := \max_{-1\le\eta\le1} |w(-L,\eta)|\to 0\qquad\text{as}~L\to\infty.
	\end{equation*}
	Taking into account \eqref{comparison-first-arc} and the boundary conditions in \eqref{difference-problem}, the boundary values of $w$ on $\partial G_L$ satisfy
	\begin{equation}\label{comparison-first-boundary-bounds}
		w \ge 0 \quad\text{on } \gamma \cup \bigl([-L,0]\times\{-1,1\}\bigr),
		\qquad
		w \ge -\varepsilon_L \quad\text{on } \{-L\}\times[-1,1].
	\end{equation}
	We first claim that
	\begin{equation}\label{comparison-first-minimum-bound}
		\min_{\overline{G_L}} w \ge -\varepsilon_L.
	\end{equation}
	Indeed, if \eqref{comparison-first-minimum-bound} failed, $w$ would attain its negative minimum over $\overline{G_L}$ at some point $Z_L \in \overline{G_L}$ with $w(Z_L) < -\varepsilon_L \le 0$. If $Z_L \in G_L$, the strong maximum principle would imply that $w$ is identically constant in $G_L$, which by continuity contradicts the nonnegativity $w \ge 0$ along $\gamma$. Hence $Z_L \in \partial G_L$. In view of \eqref{comparison-first-boundary-bounds}, $Z_L$ can only lie on $(0,r)\times\{-1\} \subset \mathcal{T}_2$, which leads to a contradiction with Lemma~\ref{strip-boundary-extrema}. This confirms \eqref{comparison-first-minimum-bound}.
	
	On the other hand, since $r > 0$, the lower detachment point $\widehat B = (0,-1)$ is strictly separated from $\gamma$, that is, $\operatorname{dist}(\widehat B, \gamma) > 0$. Taking $\rho \in \left(0, \frac{1}{2}\min\{\operatorname{dist}(\widehat B,\gamma), L_0, r, 1\}\right)$, the upper half-ball $B_\rho(\widehat B) \cap S$ is entirely contained in $G_L$ for all $L \ge L_0$. Hence, for any $0 < \eta_- < \min\{\rho, \delta\}$, the point $Z_- := (0, -1+\eta_-) \in G_L$ is independent of $L$, and the   sign condition \eqref{comparison-endpoint-signs} ensures that $w(Z_-) < 0$. Combining this with \eqref{comparison-first-minimum-bound} and sending $L \to \infty$ yields
	\[
	0 > w(Z_-) \ge \min_{\overline{G_L}} w \ge -\varepsilon_L \to 0 \qquad\text{as } L \to \infty,
	\]
	which is a contradiction. This rules out the possibility that $\overline{U_+} \cap \mathcal{T}_2 \cap \{w > 0\} \neq \varnothing$. The symmetric statement for $U_-$ is established identically, and the proof is complete.
\end{proof}

Next, we demonstrate that both components $U_+$ and $U_-$ extend infinitely downstream. More generally, we establish this property for all positive superlevel sets and negative sublevel sets.

\begin{lemma}\label{strip-level-components}
	Fix $t \ge 0$. Every non-empty connected component $V_t^+$ of $U_+ \cap \{Z\in S : w(Z) > t\}$ and every non-empty connected component $V_t^-$ of $U_- \cap \{Z\in S : w(Z) < -t\}$ satisfy
	\begin{equation*}
		\sup\{\xi : (\xi,\eta)\in V_t^+\} = \sup\{\xi : (\xi,\eta)\in V_t^-\} = +\infty.
	\end{equation*}
	In particular, by setting $t = 0$, we have
	\begin{equation*}
		\sup\{\xi : (\xi,\eta)\in U_+\} = \sup\{\xi : (\xi,\eta)\in U_-\} = +\infty.
	\end{equation*}
\end{lemma}
\begin{proof}
	We only prove the assertion for $V_t^+$, while the result for $V_t^-$ is analogous. Assume by contradiction that $\sup\{\xi : (\xi,\eta)\in V_t^+\} < \infty$. We then choose a constant $R_0>\sup\{\xi : (\xi,\eta)\in V_t^+\}$.
	Take a point $Z_0=(\xi_0,\eta_0) \in V_t^+$. Since $w(Z_0) > t $, the uniform upstream decay of $w$ in \eqref{difference-asymptotics} enables us to choose $L > 0$ such that
	\begin{equation}\label{d13}
		\xi_0> -L, \qquad \sup_{\xi \le -L, \, -1 \le \eta \le 1} |w(\xi,\eta)| < \frac{w(Z_0)}{2}.
	\end{equation}
	On the non-empty compact set $K = \overline{V_t^+} \cap ([-L,R_0]\times[-1,1])$, the continuous function $w$ attains a maximum at some point $Z_*=(\xi_*,\eta_*) \in K$. We denote this maximum value by
	\begin{equation*}
		M = w(Z_*) = \max_K w = \sup_{V_t^+} w \ge w(Z_0) > t.
	\end{equation*}
	The identity $\max_K w = \sup_{V_t^+} w$ is valid because $V_t^+ \subset \{\xi < R_0\}$, and estimate \eqref{d13} ensures that $w$ is strictly smaller than $w(Z_0)$ throughout the upstream region $\{\xi < -L\}$.
	
	We now inspect the location of $Z_*$ to derive a contradiction.
	
	First, estimate \eqref{d13} directly excludes the possibility that $Z_* \in \{\xi = -L\}$. Moreover, by the continuity of $w$, we have $\partial V_t^+ \cap S \subset \{w = t\}$, which means that $Z_*$ cannot lie on the relative boundary $\partial V_t^+ \cap S$ because $M > t$.
	
	Next, suppose that $Z_* \in V_t^+$. Since $w$ attains its maximum on the open connected set $V_t^+$, the strong maximum principle implies that $w \equiv M$ throughout $V_t^+$. This identity yields $\partial V_t^+ \cap S = \varnothing$, which in turn forces $V_t^+ = S$. This contradicts the assumption that $V_t^+ \subset \{\xi < R_0\}$.
	
	Therefore, the maximum point $Z_*$ must belong to the boundary $\partial S$. If $\xi_{*} \le 0$, the homogeneous Dirichlet condition on the nozzle wall gives $w(Z_*) = 0 < M$.
If $Z_* \in \mathcal{T}_1$, since $M > t$ and $Z_* \in \overline{V_t^+}$, by the continuity of $w$ there exists $r > 0$ such that $B_r(Z_*) \cap S \subset V_t^+$ and $w \le M$ on $B_r(Z_*) \cap \overline S$. Thus $Z_*$ is a local maximum of $w$ relative to $\overline S$ on $\mathcal{T}_1$, which contradicts Lemma~\ref{strip-boundary-extrema}.
	Finally, if $Z_* \in \mathcal{T}_2$, we would have $Z_* \in \overline{U_+} \cap \{Z\in\mathcal{T}_2 : w(Z) > 0\}$, which contradicts Lemma~\ref{strip-components}.
	
	 Since all possible locations for $Z_*$ are excluded, the proof is complete.
\end{proof}

\subsection{Proof of the uniqueness of the speed \texorpdfstring{$\lambda$}{lambda}}\label{s4.2}

To establish the uniqueness of the speed $\lambda$,  we will show that any uniform affine asymptotic formula for $w$ as $\xi\to+\infty$ forces the linear coefficient in $\xi$ to vanish and the coefficient of $\eta$ to be non-negative.

\begin{proposition}\label{strip-comparison}
	Let $w$ satisfy the sign condition \eqref{comparison-endpoint-signs}. If $A_0,B_0,D_0\in\mathbb R$ satisfy
	\begin{equation}\label{comparison-downstream}
		\lim_{\xi\to+\infty}\sup_{-1\le\eta\le1}
		|w(\xi,\eta)-A_0\xi-B_0-D_0\eta|=0,
	\end{equation}
	then
	\begin{equation}\label{comparison-coefficient-conclusion}
		A_0=0 \qquad\text{and}\qquad D_0\ge0.
	\end{equation}
\end{proposition}

\begin{proof}
	This proof is split into two steps.
	
	\textbf{Step 1}. We prove that $A_0=0$.

	If $A_0>0$,  the convergence in \eqref{comparison-downstream} implies that for any large $\xi$, and $\eta\in[-1,1]$, there holds
	\[
	w(\xi,\eta) \ge A_0\xi-|B_0|-|D_0|-1 = A_0\xi-K_0 > 0.
	\]
  It follows that $U_-$ is bounded from the right, contrary to Lemma \ref{strip-level-components}. The case $A_0<0$ is treated analogously to show that $U_+$ is bounded from the right, which again contradicts  Lemma \ref{strip-level-components}. Therefore, we obtain $A_0=0$.
	
	\textbf{Step 2}. We prove that $D_0\ge0$.
	
	Suppose, to the contrary, that $D_0<0$, and define the affine function $\ell(\eta)=B_0+D_0\eta$. By Step~1, the downstream convergence \eqref{comparison-downstream} becomes
	\begin{equation}
	\lim_{\xi\to+\infty}\sup_{-1\le\eta\le1}
	|w(\xi,\eta)-\ell(\eta)|=0.\label{d-13}
	\end{equation}
	
	Let $\delta$ be the constant in \eqref{comparison-endpoint-signs}. Choose $\eta_0\in(1-\delta,1)$ and set $P_0=(0,\eta_0)$. The upper sign condition in \eqref{comparison-endpoint-signs} gives $w(P_0)>0$. Fix $t\in(0,w(P_0)/2)$, and let $V_t$ be the connected component of $U_+\cap\{w>t\}$ containing $P_0$. Then Lemma \ref{strip-level-components} ensures that
	\begin{equation}\label{comparison-Vt-unbounded}
		\sup\{\xi:(\xi,\eta)\in V_t\}=+\infty.
	\end{equation}
	By \eqref{d-13}, we can choose $R>1$ sufficiently large such that
	\begin{equation}\label{comparison-right-error}
		\sup_{-1\le\eta\le1}|w(R,\eta)-\ell(\eta)|<\frac t3.
	\end{equation}
	By \eqref{comparison-Vt-unbounded}, $V_t$ contains a point whose first coordinate is larger than $R$. Since $V_t$ is a connected open subset of $\mathbb R^2$, it is path-connected. A path in $V_t$ from $P_0$ to such a point must intersect $\{\xi=R\}$. Hence there exists
	\begin{equation}\label{comparison-right-point}
		P_R=(R,\eta_R)\in V_t \qquad\text{with } -1<\eta_R<1.
	\end{equation}
	
	Since $D_0<0$, the function $\ell$ is strictly decreasing. Therefore, using equations \eqref{comparison-right-error} and \eqref{comparison-right-point}, we have
	\[
	w(R,\eta) > \ell(\eta)-\frac t3 \ge \ell(\eta_R)-\frac t3 >w(R,\eta_R)-\frac{2t}{3}> \frac t3
	\qquad\text{for all } -1\le\eta\le\eta_R,
	\]
	which implies that $\{R\}\times(-1,\eta_R]\subset U_+$. Consequently, $(R,-1)\in\overline{U_+}$ with $w(R,-1)\ge t/3>0$. This contradicts  Lemma~\ref{strip-components}. Therefore, we have $D_0\ge0$, and \eqref{comparison-coefficient-conclusion} follows.
\end{proof}

Finally, we will complete the proof of Theorem \ref{prop-speed} as follows.

\begin{proof}[Proof of Theorem~\ref{prop-speed}]
	The proof is divided into two steps.
	
	\textbf{Step 1}. We establish the sign condition \eqref{comparison-endpoint-signs}.
	
	Since the two jets have the same nozzle walls, we have
	\[
	h(0,1)=\tilde h(0,1)=g_1(0), \qquad
	h(0,-1)=\tilde h(0,-1)=g_2(0).
	\]
	Thus $w(0,1)=w(0,-1)=0$. Using formula \eqref{inverse-endpoint-derivatives}, we have
	$$
			w_\eta(0,1)  = \sqrt{1+g_1'(0)^2}\left(\frac1\lambda-\frac1{\tilde\lambda}\right) < 0,\quad
			w_\eta(0,-1)  = \sqrt{1+g_2'(0)^2}\left(\frac1\lambda-\frac1{\tilde\lambda}\right) < 0.
	$$
	Therefore, by the $C^1$-regularity of $w$ in Proposition~\ref{prop-inverse}, there is a constant $\delta\in(0,2/3)$ such that
	\[
	w_\eta(0,\eta)<0 \quad\text{for all } 1-\delta\le\eta\le1,
	\qquad
	w_\eta(0,\eta)<0 \quad\text{for all } -1\le\eta\le-1+\delta.
	\]
	For $1-\delta<\eta<1$, integrating from $\eta$ to $1$ gives
	\[
	w(0,\eta) = w(0,1)-\int_\eta^1w_\eta(0,\tau)\,d\tau > 0.
	\]
	For $-1<\eta<-1+\delta$, integrating from $-1$ to $\eta$ gives
	\[
	w(0,\eta) = w(0,-1)+\int_{-1}^{\eta}w_\eta(0,\tau)\,d\tau < 0.
	\]
	Hence \eqref{comparison-endpoint-signs} holds.
	
	\textbf{Step 2}. We derive a contradiction from the downstream coefficient $D_0$.
	
	Since Step 1 gives the hypothesis in  Proposition~\ref{strip-comparison}, it holds
	$$
	A_0=\tan\theta -\tan\tilde\theta=0\qquad \mathrm{and}\qquad D_0=\frac{\sec\theta}{\lambda}-\frac{\sec\tilde\theta}{\tilde{\lambda }}\ge 0.
	$$
	 Consequently, we obtain $\lambda\leq\tilde{\lambda}$, which contradicts  the  assumption $\lambda>\tilde\lambda$. Since any unequal pair can be relabeled to satisfy this inequality, it follows that $\lambda=\tilde\lambda$, which completes the proof.
\end{proof}

\section{The uniqueness of solutions to two-dimensional asymmetric incompressible jet flows}
\label{sec-main-proof}

By Theorem~\ref{prop-speed}, we have $\lambda=\widetilde\lambda$. Consequently, Proposition~\ref{prop-height-difference} shows that the height-function difference $w=h-\widetilde h$ satisfies
\begin{equation}\label{common-speed-problem}
	\left\{
	\begin{aligned}
		&\operatorname{div}(\mathcal A\nabla w)=0
		&&\text{in }S,\\
		&w(\xi,\pm1)=0
		&&\text{for }\xi\le0,\\
		 &\mathcal A\nabla w\cdot\nu_1
            =0 \ \text{and}\   w=f_1(\xi)-\tilde{f}_1(\xi)
            &&\text{on }\mathcal{T}_1,\\
            &\mathcal A\nabla w\cdot\nu_2
            =0 \ \text{and} \ w=f_2(\xi)-\tilde{f}_2(\xi)
            &&\text{on }\mathcal{T}_2.
	\end{aligned}
	\right.
\end{equation}
where $\mathcal A$ is defined by \eqref{difference-matrix} and $\nu$ denotes the outward unit normal to $S$. The purpose of this section is to prove that $w\equiv0$ on $\overline S$, and hence that the height functions are identical. Inverting the partial hodograph transformation will then give the coincidence of the physical jets and complete the proof of Theorem~\ref{thm-main}.

First, we establish an extremum principle for $w$. 

\begin{lemma}\label{lem-extrema}
	If $w$ attains a nonpositive global minimum or a nonnegative
	global maximum at a finite point of $\overline S$, then
	$w\equiv0$ in $\overline S$.
\end{lemma}

\begin{proof}
	We prove the assertion for a nonpositive global minimum.
	The assertion for a nonnegative global maximum follows by
	interchanging $h$ and $\widetilde h$.
	
	Suppose, for contradiction, that $w\not\equiv0$ in $S$ and
	that $w$ attains a nonpositive global minimum at a point
	$Z_*\in\overline S$. Since $w=0$ on the fixed walls, $w$
	cannot be a nonzero constant. The strong maximum principle
	therefore gives $w>w(Z_*)$ in $S$. In particular,
	$Z_*\notin S$. We consider the following two cases.
	
	\medskip
	\noindent\textbf{Case 1.}
	$Z_*\in\mathcal T_1\cup\mathcal T_2$.
	
	Since $\mathcal A$ is symmetric, locally H\"older continuous
	and uniformly elliptic near $Z_*$, the Hopf boundary point
	lemma \cite[Theorem~1.1]{Sabina2015} gives
	$\partial_\nu w(Z_*)<0$, where $\nu$ is the outward unit
	normal to $S$ at $Z_*$. Moreover, $w_\xi(Z_*)=0$, because
	$Z_*$ is a minimum point on the horizontal boundary. Hence
	\[
	\mathcal A(Z_*)\nabla w(Z_*)\cdot\nu
	=
	\bigl(\nu^T\mathcal A(Z_*)\nu\bigr)
	\partial_\nu w(Z_*)<0,
	\]
	which contradicts the homogeneous conormal condition in
	\eqref{common-speed-problem}.
	
	\medskip
	\noindent\textbf{Case 2.}
	$Z_*$ lies on a fixed wall or at one of the endpoints
	$\widehat A,\widehat B$.
	
	Since $w=0$ on the fixed walls and at the endpoints,
	$w(Z_*)=0$ and $w\ge0$ on $\overline S$. The strong maximum
	principle and the argument in Case~1 imply that
	$w>0$ in $S\cup\mathcal T_1\cup\mathcal T_2$. Moreover,
	$w(\widehat A)=w(\widehat B)=0$, so both endpoints are
	minimum points. We will use $\widehat A$ to derive a
	contradiction.
	
	By Proposition~\ref{prop-inverse} and
	\eqref{inverse-endpoint-derivatives}, we may choose
	$r\in(0,1)$ sufficiently small such that
	$h,\widetilde h\in C^1([-r,r]\times[1-2r,1])$ and
	$c\le h_\eta,\widetilde h_\eta\le C$ on this rectangle, for
	some constants $c,C>0$. Set
	$D_r=(-r,r)\times(1-r,1)$.
	
	We first show that there exists $c_0>0$ such that
	\begin{equation}\label{extrema-boundary-lower}
		w(\xi,\eta)\ge c_0(1-\eta)
		\qquad\text{on }\partial D_r\cap\{\eta<1\}.
	\end{equation}
	Indeed, the $C^{1,\beta}$ regularity of the fixed wall near
	$(-r,g_1(-r))$, the boundary regularity of the harmonic
	functions $u,\widetilde u$, and the inverse function theorem
	show that $\mathcal A$ is locally H\"older continuous and
	uniformly elliptic near $(-r,1)$. Since $w(-r,1)=0$ and
	$w>0$ in $S$, the Hopf boundary point lemma
	\cite[Theorem~1.1]{Sabina2015} gives $w_\eta(-r,1)<0$.
	Consequently, $w(-r,\eta)\ge c_1(1-\eta)$ for
	$1-r\le\eta<1$ and some $c_1>0$. On the bottom and right
	sides of $D_r$, the continuity and strict positivity of $w$
	give a positive lower bound. Thus
	\eqref{extrema-boundary-lower} follows.
	
	Choose $0<\varepsilon<\min\{1,c_0/C\}$ and define
	\[
	h_\varepsilon(\xi,\eta)
	=
	h\bigl(\xi,\eta-\varepsilon(1-\eta)\bigr)
	\qquad\text{on }\overline{D_r}.
	\]
   Then a direct computation gives that 
	  $h_\varepsilon$ satisfies the same height equation 
	  \eqref{equation-height}  as
	$h$ in $D_r$. 
	Moreover, $0\le h-h_\varepsilon
	\le C\varepsilon(1-\eta)$ on $\overline{D_r}$.
	Together with \eqref{extrema-boundary-lower}, this gives
	$h_\varepsilon\ge\widetilde h$ on
	$\partial D_r\cap\{\eta<1\}$. The same inequality holds on
	the top side of $D_r$, since $h_\varepsilon=h$ there and
	$w=h-\widetilde h\ge0$.
	
	  Subtracting
	the height equations for $h_\varepsilon$ and $\widetilde h$
	yields
	\[
	\operatorname{div}
	(\mathcal A_\varepsilon\nabla (h_\varepsilon-\widetilde h))=0
	\quad\text{in }D_r,\qquad
	\mathcal A_\varepsilon
	=
	\int_0^1
	D^2F\bigl((1-t)\nabla\widetilde h
	+t\nabla h_\varepsilon\bigr)\,dt.
	\]
	The local $C^1$ bounds and the positivity of
	$(h_\varepsilon)_\eta,\widetilde h_\eta$ imply that
	$\mathcal A_\varepsilon$ is bounded and uniformly elliptic
	on $\overline{D_r}$. The maximum principle therefore gives
	$h_\varepsilon\ge\widetilde h$ in $\overline{D_r}$.
	
	Finally, $h_\varepsilon(\widehat A)=\widetilde h(\widehat A)$.
	Using \eqref{inverse-endpoint-derivatives} together with
	$\lambda=\widetilde\lambda$, we obtain
	\[
	0\le
	\widetilde h_\eta(\widehat A)
	-(1+\varepsilon)h_\eta(\widehat A)
	=-\varepsilon h_\eta(\widehat A)<0,
	\]
	a contradiction. The proof is complete.
\end{proof}

With the aid of the maximum principle for $w$ established in
Lemma \ref{lem-extrema} and the asymptotic behavior of
$w$ given in Proposition \ref{prop-height-difference}, we now prove the
uniqueness of the height functions.

\begin{proposition}\label{prop-fixed-speed-height}
	If $\lambda=\widetilde\lambda$, then
	\[
	h=\widetilde h\quad\text{in }\overline S,\qquad
	\tan\theta=\tan\widetilde\theta,\qquad
	b_i=\widetilde b_i\quad\text{for }i=1,2.
	\]
\end{proposition}

\begin{proof}
	We first show that $A_0=0$. If not, without loss of generality, we suppose that $A_0>0$. By the uniform downstream
	expansion in \eqref{difference-asymptotics},
	\begin{equation}\label{e15}
		\inf_{-1\le\eta\le1}w(\xi,\eta)
		\geq \frac{A_0\xi+B_0-|D_0|}2 \to+\infty
		\qquad\text{as }\xi\to+\infty.
	\end{equation}
	By the upstream limit in \eqref{difference-asymptotics} and
	the downstream limit in \eqref{e15}, $w$ is bounded from
	below in $S$. Since $w(\widehat A)=0$, we have
	$\inf_S w\le0$. If $\inf_S w=0$, then $w$ attains its
	global minimum at $\widehat A$. If $\inf_S w<0$, the two
	asymptotic limits imply that $w$ attains its global minimum
	in a finite closed substrip of $\overline S$. In either
	case, Lemma~\ref{lem-extrema} gives $w\equiv0$, which
	contradicts \eqref{e15}.

	Therefore, $A_0=0$, and
	\eqref{difference-asymptotic-coefficients} gives
	$\tan\theta=\tan\widetilde\theta$. Since
	$\theta,\widetilde\theta\in(-\pi/2,\pi/2)$, we have
	$\theta=\widetilde\theta$.
	Together with $\lambda=\widetilde\lambda$, the same
	formula \eqref{difference-asymptotic-coefficients} yields $D_0=0$.
	
	With the aid of $A_0=D_0=0$, the asymptotic behavior \eqref{difference-asymptotics}
	reduces to
	\[
	\lim_{\xi\to+\infty}\sup_{-1\le\eta\le1}
	|w(\xi,\eta)-B_0|=0.
	\]
 Using Lemma  \ref{lem-extrema} and the asymptotic behaviors of $w$, we can argue as before to show that $B_0=0$.
	Consequently,
	\eqref{difference-asymptotic-coefficients} gives
	\[
	b_1+b_2=\widetilde b_1+\widetilde b_2.
	\]
	Moreover, \eqref{jet-flux-width}, together with
	$\tan\theta=\tan\widetilde\theta$ and
	$\lambda=\widetilde\lambda$, implies
	\[
	b_1-b_2=\widetilde b_1-\widetilde b_2.
	\]
	Hence, $b_i=\widetilde b_i$ for $i=1,2$.
	
	Since $A_0=B_0=D_0=0$,
	\eqref{difference-asymptotics} now gives
	\begin{equation}\label{zero-both-ends}
		\lim_{\xi\to\pm\infty}\sup_{-1\le\eta\le1}
		|w(\xi,\eta)|=0.
	\end{equation}
	It follows from  Lemma  \ref{lem-extrema} that $w\equiv0$, and hence
	$h=\widetilde h$ on $\overline S$.
\end{proof}

Having established the identity of the height functions on the fixed strip, we are now ready to translate this equality back to the original physical coordinates and complete the proof of Theorem~\ref{thm-main}.

\begin{proof}[Completion of the proof of Theorem~\ref{thm-main}]
	By Theorem~\ref{prop-speed}, the two jets have the identical free-boundary flow speed $\lambda = \widetilde\lambda$. In turn, Proposition~\ref{prop-fixed-speed-height} ensures that $h \equiv \widetilde h$ on $\overline S$ and that the downstream asymptotic constants coincide, namely $\tan\theta = \tan\widetilde\theta$ and $b_i = \widetilde b_i$ for $i=1,2$. Since $\theta,\widetilde\theta\in(-\pi/2,\pi/2)$, this also gives $\theta=\widetilde\theta$.
	
	We first show that the fluid domains and the free boundaries coincide. Recall that the fluid regions of the two jets are given by
	\[
	\Omega = \{(x,y)\in\mathbb R^2: a(x) < y < b(x)\},
	\qquad
	\widetilde\Omega = \{(x,y)\in\mathbb R^2: \widetilde a(x) < y < \widetilde b(x)\},
	\]
	where the boundary functions satisfy
	\[
	a(x) = h(x,-1), \qquad b(x) = h(x,1),
	\]
	and
	\[
	\widetilde a(x) = \widetilde h(x,-1), \qquad \widetilde b(x) = \widetilde h(x,1),
	\]
	for all $x \in \mathbb R$. Since $h \equiv \widetilde h$ on $\overline S$, we immediately obtain $a(x) = \widetilde a(x)$ and $b(x) = \widetilde b(x)$ for all $x \in \mathbb R$. In particular, the two fluid domains are identical, $\Omega = \widetilde\Omega$. Restricting this equality to $x \ge 0$ yields $f_i(x) = \widetilde f_i(x)$ for $i=1,2$, so that the upper and lower free boundaries coincide as well.
	
	Next, we establish the identity of the stream functions. For any point $(x,y) \in \Omega$, let $\eta = u(x,y) \in (-1,1)$. By the definition of the height function, we have $y = h(x,\eta)$. Since $h = \widetilde h$, we obtain $y = \widetilde h(x,\eta)$, which in turn implies that
	\[
	\widetilde u(x,y) = \widetilde u(x, \widetilde h(x,\eta)) = \eta = u(x,y).
	\]
	Therefore, $u = \widetilde u$ in $\Omega$, and by the continuity of both stream functions up to the boundary, this equality extends to the closure $\overline\Omega$.
	
	Finally, we recover the velocity and pressure fields. In terms of the normalized variables introduced in Subsection~\ref{sec-setting}, the velocity fields of the two flows are given by
	\[
	(u_1, u_2) = (u_y, -u_x),
	\qquad
	(\widetilde u_1, \widetilde u_2) = (\widetilde u_y, -\widetilde u_x).
	\]
	Since $u \equiv \widetilde u$ in $\Omega$, their spatial gradients are identical, which immediately gives $(u_1, u_2) \equiv (\widetilde u_1, \widetilde u_2)$ in $\Omega$. Furthermore, the fluid pressures satisfy
	\[
	p = p_{\mathrm{atm}} + \frac{\lambda^2 - |\nabla u|^2}{2},
	\qquad
	\widetilde p = p_{\mathrm{atm}} + \frac{\widetilde\lambda^2 - |\nabla\widetilde u|^2}{2}.
	\]
	Because $\lambda = \widetilde\lambda$ and $\nabla u = \nabla\widetilde u$, the two pressures also coincide throughout the common fluid domain. This completes the proof of Theorem~\ref{thm-main}.
\end{proof}

\appendix

\makeatletter
\renewcommand{\@seccntformat}[1]{}
\makeatother

\section{Appendix. A maximum principle on unbounded domains}
\label{app-max}
This appendix establishes the comparison principle used for the vertical
translates of $u$ in Lemma~\ref{lem-monotonicity} and for the
far-field estimate in Lemma~\ref{lem-farfield}.

\begin{lemma}\label{lem-unbounded}
	Let $D \subset \mathbb{R}^2$ be an open set, and suppose that
	$\overline{B_r(Z_0)} \subset \mathbb{R}^2 \setminus D$ for some
	$Z_0 \in \mathbb{R}^2$ and $r > 0$. If $v \in C(\overline{D}) \cap C^2(D)$ satisfies
	\[
	\Delta v = 0 \quad \text{in } D,
	\qquad
	\sup_D v < \infty,
	\qquad
	v \le 0 \quad \text{on } \partial D,
	\]
	then $v \le 0$ in $D$.
\end{lemma}

\begin{proof}
	Set $M = \max\{0, \sup_D v\}$. Fix an arbitrary point $X_0 \in D$ and $\varepsilon > 0$. We choose $R > 0$ sufficiently large such that
	\[
	R > |X_0 - Z_0|>r \quad \text{and} \quad \varepsilon \log(R/r) \ge M.
	\]
	Let $D_R$ be the connected component of $D \cap B_R(Z_0)$ containing $X_0$, and consider the auxiliary function
	\[
	v_\varepsilon(X) = v(X) - \varepsilon \log\left(\frac{|X - Z_0|}{r}\right).
	\]
	Since $Z_0 \notin \overline{D}$, the function $X \mapsto \log|X - Z_0|$ is smooth and harmonic in $D$, which implies that $\Delta v_\varepsilon = 0$ in $D_R$.
	
	To examine the boundary values, note that $\partial D_R \subset (\partial D \cap \overline{B_R(Z_0)}) \cup (\overline{D} \cap \partial B_R(Z_0))$. On $\partial D \cap \overline{B_R(Z_0)}$, the condition $\overline{B_r(Z_0)} \subset \mathbb{R}^2 \setminus D$ implies that $|X - Z_0| \ge r$. Hence, combining this with $v \le 0$ on $\partial D$ yields
	\[
	v_\varepsilon(X) \le -\varepsilon \log\left(\frac{|X - Z_0|}{r}\right) \le 0.
	\]
	On the remaining portion $\overline{D} \cap \partial B_R(Z_0)$, we have $|X - Z_0| = R$, whence
	\[
	v_\varepsilon(X) = v(X) - \varepsilon \log(R/r) \le M - \varepsilon \log(R/r) \le 0.
	\]
	Consequently, $v_\varepsilon \le 0$ on $\partial D_R$. The classical maximum principle on bounded domains \cite[Theorem~3.1]{GT} then gives $v_\varepsilon \le 0$ in $D_R$. In particular, since $X_0 \in D_R$,
	\[
	v(X_0) \le \varepsilon \log\left(\frac{|X_0 - Z_0|}{r}\right).
	\]
	Taking the limit as $\varepsilon \to 0$ yields $v(X_0) \le 0$. Since $X_0 \in D$ was chosen arbitrarily, the proof is complete.
\end{proof}

\section*{Declarations}

\noindent\textbf{Conflict of interest.}
The authors declare that they have no conflict of interest.

\noindent\textbf{Data availability.}
No data were used or generated in this study.

\begingroup
\setlength{\parskip}{0pt}

\endgroup

\begin{thebibliography}{99}
	\addcontentsline{toc}{section}{References}
	\setlength{\itemsep}{2pt}
	\setlength{\parsep}{0pt}
	
	\bibitem{AC81}
	H.~W.~Alt and L.~A.~Caffarelli,
	Existence and regularity for a minimum problem with free boundary,
	\emph{J. Reine Angew. Math.} \textbf{325} (1981), 105--144.
	
	\bibitem{ACF}
	H.~W.~Alt, L.~A.~Caffarelli, and A.~Friedman,
	Asymmetric jet flows,
	\emph{Comm. Pure Appl. Math.} \textbf{35} (1982), no.~1, 29--68.
	
	\bibitem{ACFgravity}
	H.~W.~Alt, L.~A.~Caffarelli, and A.~Friedman,
	Jet flows with gravity,
	\emph{J. Reine Angew. Math.} \textbf{331} (1982), 58--103.
	
	\bibitem{ACFaxis}
	H.~W.~Alt, L.~A.~Caffarelli, and A.~Friedman,
	Axially symmetric jet flows,
	\emph{Arch. Ration. Mech. Anal.} \textbf{81} (1983), 97--149.
	
	\bibitem{ACFtwofluidsI}
	H.~W.~Alt, L.~A.~Caffarelli, and A.~Friedman,
	Jets with two fluids. I. One free boundary,
	\emph{Indiana Univ. Math. J.} \textbf{33} (1984), 213--247.
	
	\bibitem{ACFtwofluidsII}
	H.~W.~Alt, L.~A.~Caffarelli, and A.~Friedman,
	Jets with two fluids. II. Two free boundaries,
	\emph{Indiana Univ. Math. J.} \textbf{33} (1984), 367--391.
	
	\bibitem{ACFcompressible}
	H.~W.~Alt, L.~A.~Caffarelli, and A.~Friedman,
	Compressible flows of jets and cavities,
	\emph{J. Differential Equations} \textbf{56} (1985), no.~1, 82--141.
	
	\bibitem{Berg1962}
	P.~W.~Berg,
	The existence of subsonic Helmholtz flows of a compressible fluid,
	\emph{Comm. Pure Appl. Math.} \textbf{15} (1962), 289--347.
	
	\bibitem{BZ1957}
	G.~Birkhoff and E.~H.~Zarantonello,
	\emph{Jets, Wakes, and Cavities},
	Academic Press, New York, 1957.
	
	\bibitem{CDW}
	J.~Cheng, L.~Du, and Y.~Wang,
	The uniqueness of the asymmetric jet flow,
	\emph{J. Differential Equations} \textbf{269} (2020), no.~4, 3794--3815.
	
	\bibitem{CDXcompressible19}
	J.~Cheng, L.~Du, and W.~Xiang,
	Compressible subsonic jet flows issuing from a nozzle of arbitrary cross-section,
	\emph{J. Differential Equations} \textbf{266} (2019), no.~9, 5318--5359.
	
	\bibitem{CDX19}
	J.~Cheng, L.~Du, and W.~Xiang,
	Incompressible jet flows in a de Laval nozzle with smooth detachment,
	\emph{Arch. Ration. Mech. Anal.} \textbf{232} (2019), 1031--1072.
	
	
	\bibitem{CGX25}
	J.~Cheng, Z.~Gu, and W.~Xiang,
	Axially symmetric incompressible jet flows with vorticity,
	\emph{J. Lond. Math. Soc.} \textbf{112} (2025), no.~6, Paper No.~e70405.
	
	\bibitem{ChengHuang23}
	J.~Cheng and J.~Huang,
	Incompressible jet flow past an obstacle,
	\emph{SIAM J. Math. Anal.} \textbf{55} (2023), no.~6, 7077--7112.
	
	\bibitem{ChengPei24}
	J.~Cheng and Z.~Pei,
	The well-posedness of incompressible rotational jet flows with gravity,
	\emph{Calc. Var. Partial Differential Equations} \textbf{63} (2024), Paper No.~50.
	
	
	
	\bibitem{DuTang21}
	L.~Du and Y.~Tang,
	The steady incompressible jet flows issuing from a finitely long nozzle,
	\emph{J. Differential Equations} \textbf{280} (2021), 355--374.
	
	
	\bibitem{Finn1956}
	R.~Finn,
	Some theorems on discontinuous plane fluid motions,
	\emph{J. Analyse Math.} \textbf{4} (1954/56), 246--291.
	
	\bibitem{Friedman1984}
	A.~Friedman,
	Free boundary problem in fluid dynamics,
	\emph{Ast\'erisque} \textbf{118} (1984), 55--67.
	
	\bibitem{Gilbarg1952}
	D.~Gilbarg,
	Uniqueness of axially symmetric flows with free boundaries,
	\emph{J. Rational Mech. Anal.} \textbf{1} (1952), 309--320.
	
	\bibitem{Gilbarg1960}
	D.~Gilbarg,
	Jets and cavities,
	in \emph{Handbuch der Physik}, vol.~IX, Fluid Dynamics III,
	Springer-Verlag, Berlin, 1960, pp.~311--445.
	
	\bibitem{GT}
	D.~Gilbarg and N.~S.~Trudinger,
	\emph{Elliptic Partial Differential Equations of Second Order},
	Classics in Mathematics, Springer-Verlag, Berlin, 2001.
	Reprint of the 1998 edition.
	
	\bibitem{Kirchhoff1869}
	G.~Kirchhoff,
	Zur Theorie freier Fl\"ussigkeitsstrahlen,
	\emph{J. Reine Angew. Math.} \textbf{70} (1869), 289--298.
	
	\bibitem{KL23}
	V.~Kozlov and E.~Lokharu,
	Global bifurcation and highest waves on water of finite depth,
	\emph{Arch. Ration. Mech. Anal.} \textbf{247} (2023), Paper No.~98.
	
	\bibitem{Leray1935}
	J.~Leray,
	Les probl\`emes de repr\'esentation conforme d'Helmholtz; th\'eories des sillages et des proues, I, II,
	\emph{Comment. Math. Helv.} \textbf{8} (1935), 149--180 and 250--263.
	
	\bibitem{LSTX25}
	Y.~Li, W.~Shi, L.~Tang, and C.~Xie,
	Variational structure and two-dimensional subsonic jet flows for compressible Euler system with general incoming flows,
	\emph{Arch. Ration. Mech. Anal.} \textbf{249} (2025), Paper No.~49.
	
	\bibitem{Sabina2015}
	J.~C.~Sabina de Lis,
	Hopf maximum principle revisited,
	\emph{Electron. J. Differential Equations}
	(2015), no.~115, 1--9.
	
	\bibitem{Serrin1952Comparison}
	J.~Serrin,
	Two hydrodynamical comparison theorems,
	\emph{J. Rational Mech. Anal.} \textbf{1} (1952), 563--572.
	
	\bibitem{Serrin1952Uniqueness}
	J.~Serrin,
	Uniqueness theorems for two free boundary problems,
	\emph{Amer. J. Math.} \textbf{74} (1952), 492--506.
	
	\bibitem{Serrin1953}
	J.~Serrin,
	On plane and axially symmetric free boundary problems,
	\emph{J. Rational Mech. Anal.} \textbf{2} (1953), 563--575.
	
	\bibitem{WangXin19}
	C.~Wang and Z.~Xin,
	On an elliptic free boundary problem and subsonic jet flows for a given surrounding pressure,
	\emph{SIAM J. Math. Anal.} \textbf{51} (2019), no.~2, 1014--1045.
	
	\bibitem{Weinstein1924}
	A.~Weinstein,
	Ein hydrodynamischer Unit\"atssatz,
	\emph{Math. Z.} \textbf{19} (1924), 265--275.
	
	\bibitem{Weinstein1930}
	A.~Weinstein,
	Zur Theorie der Fl\"ussigkeitsstrahlen,
	\emph{Math. Z.} \textbf{31} (1930), 424--433.
	
	\bibitem{Wu1972}
	T.~Y.~Wu,
	Cavity and wake flows,
	\emph{Annu. Rev. Fluid Mech.} \textbf{4} (1972), 243--284.
	
\end{thebibliography}
\end{document}